\documentclass[final]{article}
\usepackage{ifthen} 

\newboolean{preprint} 
\setboolean{preprint}{false} 

\newboolean{draft} 
\setboolean{draft}{false} 

\newboolean{final} 
\setboolean{final}{true}

\ifpreprint
\setboolean{final}{false} 
\fi
\ifdraft
\setboolean{final}{false} 
\fi

\usepackage{a4,
srcltx,
enumerate,
amsmath,amssymb}
\input{xy}
\xyoption{all}

\usepackage{amsthm,
verbatim,
varioref,
doc,
graphicx}
\usepackage[latin1]{inputenc}

\usepackage{cyrillic} 

\usepackage{srcltx,amssymb,enumerate,stackrel}
\usepackage{xr-hyper}
\usepackage{hyperref}
\makeatletter
\@addtoreset{equation}{chapter}
\@addtoreset{equation}{section}
\makeatother

		\newcommand{\category}[1]{\mathbf{#1}}
	
		\newcommand{\Com}{\category{Com}} 
    \newcommand{\K}{\category{K}} 
    \newcommand{\D}{\category {D}} 

    \newcommand{\gm}{\mathrm{gm}}

    \newcommand{\Ab}{\category {Ab}} 
    \newcommand{\Sets}{\category {Sets}} 
    
		\newcommand{\Sch}{\category{Sch}} 

    \newcommand{\Sm}{\category{Sm}}

		\newcommand{\Spt}{\category{Spt}} 

		\newcommand{\Shvv}{\category{Shv}} 

    \newcommand{\DM}{\category {DM}}
    
    \newcommand{\DMBei}{\DM_{\Beilinson}}
    \newcommand{\DMBeic}{\DM_{\Beilinson, c}}

    \newcommand{\DMgmbeff}{\DM^\eff_\gm} 
    \newcommand{\DMgmeff}{{\DMgmbeff}} 
    \newcommand{\DMtrgm}{{\DM_{\Beilinson, \mathrm{tr}}}} 
    \newcommand{\DMgm}{\DM_\gm} 
    \newcommand{\DTM}{\category{DTM}} 
    \newcommand{\MTM}{\category{MTM}} 
    \newcommand{\DATM}{\category{DATM}} 
    \newcommand{\MATM}{\category{MATM}} 

    \newcommand{\DbH}{\D^\bound_\mathrm{H}}

    \newcommand{\MM}{\category {MM}} 
    \newcommand{\PureMot}{\category {M}} 
    \newcommand{\MHS}{\category {MHS}} 

\newcommand{\DbRQdet}{\D^\bound(\underline \RR)^{\Q-\det}} 

\newcommand{\SH}{\category{SH}} 

	    \newtheoremstyle{Normal}
  {}
  {}
  {}
  {}
  {\bfseries}
  {.}
  { }
  {}

    \theoremstyle{Normal} 

    \newtheorem{Defi}{Definition}[section]

    \newtheorem{Conj}[Defi]{Conjecture}
    \newtheorem{Bem}[Defi]{Remark}
    \newtheorem{Bsp}[Defi]{Example}
     
    \newtheorem{Axio}[Defi]{Axiom}
    \newtheorem{Ques}[Defi]{Question}
    \newtheorem{Nota}[Defi]{Notation}
    \newtheorem{Comp}[Defi]{Complements}
    
		\theoremstyle{remark}
	
    \theoremstyle{plain} 

    \newtheorem{Satz}[Defi]{Proposition}
    \newtheorem{DefiTheo}[Defi]{Definition and Theorem}
    \newtheorem{Theo}[Defi]{Theorem}
    \newtheorem{Folg}[Defi]{Corollary}
    \newtheorem{Lemm}[Defi]{Lemma}
    \newtheorem{DefiLemm}[Defi]{Definition and Lemma}

\newcommand{\refit}[1]{(\ref{item_#1})}

\newcommand{\refsect}[1]{Section \ref{sect_#1}}
\newcommand{\conj}{\begin{Conj}} 			\newcommand{\xconj}{\end{Conj}}											    \newcommand{\refco}[1]{Conjecture \ref{conj_#1}}
\newcommand{\ques}{\begin{Ques}} 			\newcommand{\xques}{\end{Ques}}											    
\newcommand{\axio}{\begin{Axio}} 			\newcommand{\xaxio}{\end{Axio}}											    \newcommand{\refax}[1]{Axiom \ref{axio_#1}}
\newcommand{\bem}{\begin{Bem}} 			\newcommand{\xbem}{\end{Bem}}											    \newcommand{\refbe}[1]{Remark \ref{bem_#1}}
\newcommand{\rema}{\begin{Bem}} 			\newcommand{\xrema}{\end{Bem}}											    \newcommand{\refre}[1]{Remark \ref{rema_#1}}
\newcommand{\defi}{\begin{Defi}} 			\newcommand{\xdefi}{\end{Defi}}										\newcommand{\refde}[1]{Definition \ref{defi_#1}}
\newcommand{\defitheo}{\begin{DefiTheo}} \newcommand{\xdefitheo}{\end{DefiTheo}} 
\newcommand{\defilemm}{\begin{DefiLemm}} \newcommand{\xdefilemm}{\end{DefiLemm}} 

\newcommand{\lemm}{\begin{Lemm}}			\newcommand{\xlemm}{\end{Lemm}}											\newcommand{\refle}[1]{Lemma \ref{lemm_#1}}
\newcommand{\nota}{\begin{Nota}}			\newcommand{\xnota}{\end{Nota}}											\newcommand{\refno}[1]{Notation \ref{nota_#1}}
\newcommand{\comp}{\begin{Comp}}			\newcommand{\xcomp}{\end{Comp}}											

\newcommand{\satz}{\begin{Satz}}			\newcommand{\xsatz}{\end{Satz}}										\newcommand{\refsa}[1]{Proposition \ref{satz_#1}}

\newcommand{\prop}{\begin{Satz}}			\newcommand{\xprop}{\end{Satz}}										\newcommand{\refpr}[1]{Proposition \ref{prop_#1}}
\newcommand{\theo}{\begin{Theo}}			\newcommand{\xtheo}{\end{Theo}}											\newcommand{\refth}[1]{Theorem \ref{theo_#1}}
\newcommand{\bsp}{\begin{Bsp}}				\newcommand{\xbsp}{\end{Bsp}}												\newcommand{\refbs}[1]{Example \ref{bsp_#1}}
\newcommand{\exam}{\begin{Bsp}}				\newcommand{\xexam}{\end{Bsp}}												\newcommand{\refex}[1]{Example \ref{exam_#1}}
\newcommand{\folg}{\begin{Folg}}				\newcommand{\xfolg}{\end{Folg}}
\newcommand{\coro}{\folg}				\newcommand{\xcoro}{\xfolg} 			\newcommand{\refcor}[1]{Corollary \ref{coro_#1}}
\newcommand{\cor}{\begin{Folg}}				\newcommand{\xcor}{\end{Folg}}
\newcommand{\mycomment}{}

\newcommand{\eqnarra}{\begin{eqnarray}}				\newcommand{\xeqnarra}{\end{eqnarray}}
\newcommand{\eqnarr}{\begin{eqnarray*}}				\newcommand{\xeqnarr}{\end{eqnarray*}}

\newcommand{\eqn}{\begin{equation}} 		\newcommand{\xeqn}{\end{equation}}
\newcommand{\refeq}[1]{(\ref{eqn_#1})}

    \newcommand{\gl}{\index}

    \newcommand{\Dph}[1]{\emph{#1}}
  \newcommand{\Def}[1]{\Dph{#1}\gl{#1}}

\newcommand{\mylabel}[1]{\label{#1}}
\ifdraft
\renewcommand{\mylabel}[1]{\label{#1}\framebox[1.1\width]{$#1$}} 
\fi

\newcommand{\citefcohoAxiom}[1]{\cite[Axiom #1]{Scholbach:fcoho}}

\newcommand{\refaxcohomdim}{\citefcohoAxiom{4.1.}}

\newcommand{\refaxtstructureRealizations}{\citefcohoAxiom{4.5.}}

\newcommand{\fcohoAxiomSections}{Section 4}

\newcommand{\fcohoSectionThree}{Section 3}

\newcommand{\fcohoLemmcohomdimOF}{Lemma 5.2}
\newcommand{\fcohoSatzgeneratorsoverOF}{Prop. 5.6}

\newcommand{\fcohoSectionGenericIntermediateExtension}{Section 5.4}
\newcommand{\fcohoSectionIntermediateExtensionEll}{Section 5.5}
\newcommand{\fcohoLemmprepinjectivity}{Lemma 6.9}
\newcommand{\fcohoTheoSummaryHOneF}{Theorem 6.11}

\newcommand{\mixedLemmSplit}{Lemma 2.5}
\newcommand{\mixedTheoExactness}{Theorem 4.2}
\newcommand{\mixedSatzCohomdimAT}{Proposition 4.4}

\newcommand{\status}[1]{}

\newcommand{\OF}{{\mathcal{O}_F}}

\newcommand{\dual}{^\vee} 

\renewcommand{\log}{\mathrm{log}\,}
\newcommand{\C}{\mathcal{C}}

\newcommand{\R}{\mathrm{R}}
\newcommand{\RG}{{\R} {\Gamma}}
\newcommand{\RGw}{\RG_\mathrm{w}} 
\newcommand{\Hhat}{\widehat{\H}}
\newcommand{\Hw}{\H_\mathrm{w}} 
\newcommand{\HD}{\H_\mathrm{D}} 
\newcommand{\HB}{\mathrm{H}_{\Beilinson}}
\newcommand{\HBetti}{\mathrm{H}_{\Betti}}
\newcommand{\HBR}{\H_{\Betti, \RR}} 
\newcommand{\HBQ}{\H_{\Betti, \Q}} 
\newcommand{\HdR}{\H_\dR} 

\newcommand{\lr}{{\longrightarrow}}
\renewcommand{\r}{\rightarrow}
\renewcommand{\l}{{\leftarrow}}

\newcommand{\cC}{\mathcal{C}}
\renewcommand{\t}{{\otimes}}

\newcommand{\A}{\mathbb{A}^1}
\renewcommand{\P}[1][1]{\mathbb P^{#1}}

\newcommand{\Q}{\mathbb{Q}}

\newcommand{\Ql}{{\mathbb{Q}_{{\ell}}}}

\newcommand{\CC}{\mathbb{C}}
\newcommand{\RR}{\mathbb{R}}
\newcommand{\Z}{\mathbb{Z}}

\newcommand{\Deligne}{\mathrm{D}}
\newcommand{\Fp}{\mathbb{F}_p}

\newcommand{\pp}{\mathfrak{p}}
\newcommand{\Fpp}{{\mathbb{F}_{\pp}}}

\newcommand{\Fq}{{\mathbb{F}_q}}
\newcommand{\Fpq}{{\overline{\mathbb{F}}_p}}
\newcommand{\Fqq}{{\overline{\Fq}}}

\renewcommand{\H}{\mathrm{H}}

\newcommand{\compact}{\mathrm{c}} 
\newcommand{\Hc}{\H_\compact}

\newcommand{\RGc}{\RG_\compact}
\newcommand{\p}{{{}^\mathrm{p}}}
\newcommand{\pH}{{\p\H}}
\newcommand{\HH}{\mathbb{H}}
\newcommand{\x}{{\times}}
\newcommand{\ol}[1]{{\overline{#1}}} 
\newcommand{\one}{\mathbf{1}}
\newcommand{\onehat}{\widehat \one}

\newcommand{\Beweis}{{\normalfont} \textbf{Proof}}

\newcommand{\lc}{\textit{loc.~cit.}}

\newcommand{\lcs}{\textit{loc.~cit.}\ }
\newcommand{\ocs}{\textit{op.~cit.}\ }
\newcommand{\Spec}{\mathrm{Spec}\text{ }}
\newcommand{\SpecOF}{{\Spec} {\OF}}

\newcommand{\SpecZ}{{\Spec} {\Z}}

\newcommand{\SpecFp}{{\Spec} {\Fp}}

\renewcommand{\det}{\operatorname{det}}
\newcommand{\Fr}{\operatorname{Fr}} 

\newcommand{\SpecFpp}{{\Spec} {\Fpp}}

\newcommand{\op}{\mathrm{op}}
\newcommand{\rk}{\operatorname{rk}} 

\newcommand{\eff}{\mathrm{eff}}

\newcommand{\red}{\mathrm{red}}

\newcommand{\Gm}{\mathbb{G}_m}

\newcommand{\ord}{\operatorname{ord}} 
\newcommand{\gr}{\operatorname{gr}} 

    \def\acute{\mathaccent"7013 }
    \newcommand{\et}{\mathrm{\acute et}}

\newcommand{\an}{\mathrm{an}} 
\newcommand{\h}{\mathrm{h}}

\newcommand{\ch}{\operatorname{ch}} 
\newcommand{\dR}{\mathrm{dR}}
\newcommand{\Betti}{\mathrm{B}}
\newcommand{\alg}{\mathrm{alg}}

\newcommand{\wt}{\operatorname{wt}} 
\newcommand{\num}{\mathrm{num}} 
\renewcommand{\hom}{\mathrm{hom}} 
\newcommand{\rat}{\mathrm{rat}} 
\newcommand{\Hom}{\mathrm{Hom}}
\newcommand{\ind}{\mathrm{ind}} 
\newcommand{\IHom}{\underline{\Hom}}
\newcommand{\id}{\mathrm{id}}
\renewcommand{\gr}{\operatorname{gr}} 

\newcommand{\im}{\operatorname{im}} 

\newcommand{\M}{\operatorname{M}} 
\renewcommand{\h}{\operatorname{h}} 
\newcommand{\Mc}{\M_\mathrm c} 
\newcommand{\Mgm}{{\M}_{{\gm}}}

\newcommand{\Gal}{\mathrm{Gal}}
\newcommand{\CH}{\mathrm{CH}}
\newcommand{\CHhat}{\widehat \CH{}}
\newcommand{\Khat}{\widehat{\mathrm{K}}}

\renewcommand{\mod}[1]{\ \ (\mathrm{mod \ }#1)}
\newcommand{\CHhatGS}{\widehat{\CH}_{\mathrm{GS}{}}}

\newcommand{\CHtilde[1]}{\CHhatGS^{#1}}
\newcommand{\bound}{\mathrm{b}} 
\newcommand{\coker}{\operatorname{coker}}
\newcommand{\Ext}{\operatorname{Ext}}

\newcommand{\cone}{\operatorname{cone}}
\newcommand{\Tot}{\mathrm{Tot}}
\newcommand{\Id}{\mathrm{Id}}

\newcommand{\loc}[2]{[}

\newcommand{\pr}{\begin{proof}[\Beweis: ]}
\newcommand{\pf}{\pr}
\newcommand{\xpf}{\end{proof}}
\newcommand{\xproof}{\end{proof}}
\newcommand{\cyrrm}{\fontencoding{OT2}\selectfont\textcyrup}
\DeclareMathOperator{\Beilinson}{{\cyrrm{B}}} 
\newcommand{\twi}[1]{\{#1\}}

\newcommand{\BGL}{\mathrm{BGL}}
\newcommand{\BGLhat}{\widehat{\BGL}}

\makeindex

\numberwithin{equation}{section}

\begin{document}

\title{Special $L$-values of geometric motives}
\author{Jakob Scholbach \footnote{Universit{\"a}t M{\"u}nster, Germany}}
\maketitle

\begin{abstract}
This paper proposes a conceptual unification of Beilinson's conjecture about special $L$-values for motives over $\Q$, the Tate conjecture over $\Fp$ and Soul\'e's conjecture about pole orders of $\zeta$-functions of schemes over $\Z$. We conjecture the following: the order of $L(M, s)$ at $s=0$ is given by the negative Euler characteristic of motivic cohomology of 
$M\dual(-1)$. Up to a nonzero rational factor, the $L$-value at $s=0$ is given by the determinant of the pairing of Arakelov motivic cohomology of $M$ with the motivic homology of $M$:
$$L^*(M, 0) \equiv \prod_{i \in \Z} \det(\H_{i-2}(M, -1)  \t \Hhat^{i}(M) \r \RR)^{(-1)^{i+1}} \mod{\Q^\x}.$$ 
Under standard assumptions concerning mixed motives over $\Q$, $\Fp$, and $\Z$, this conjecture is equivalent to the conjunction of the above-mentioned conjectures of Beilinson, Tate, and Soul\'e. We use this to unconditionally prove the Beilinson conjecture for all Tate motives and, up to an $n$-th root of a rational number, for all Artin-Tate motives.
\end{abstract}

\ifpreprint
\tableofcontents
\fi

In this paper, we study special values of $L$-functions of geometric motives over $\Z$. This contains both $L$-functions over $\Q$ and Hasse-Weil $\zeta$-functions of schemes $X$ of finite type over $\Z$ (Propositions \ref{satz_corrj}, \ref{satz_HasseWeil}):
\eqn \mylabel{eqn_LQZ}
L_\Q(M_\eta, s)^{-1} = L_\Z (\eta_{!*} M_\eta[1], s), 
\xeqn
$$\zeta(X, s) = L(\Mc(X), s).$$
Here $M_\eta$ is a mixed motive over $\Q$, $\eta_{!*}$ is a generic intermediate extension functor similar to the one familiar in perverse sheaf theory, and $\Mc(X)$ denotes the motive with compact support. 

Our conjecture on special $L$-values is as follows:

\conj \mylabel{conj_summary} (see Conjectures \ref{conj_motcomp} and \ref{conj_L}) 
Let $M$ be any geometric motive over $\Z$. We conjecture that pole orders are given by the negative Euler characteristic of motivic cohomology of $M\dual(-1)$:
\eqn
\mylabel{eqn_pole.order}
\ord_{s=0} L(M, s) = -\chi (M\dual(-1)).
\xeqn
We conjecture that the \emph{Arakelov intersection pairing}, which is the natural pairing of $\RR$-vector spaces
$$\pi_M: \underbrace{\Hom (\one(-1)[-2], M)}_{=:\H_{-2}(M, -1)} \x \underbrace{\Hom (M, \onehat)}_{=:\Hhat^0(M)} \stackrel{\circ}\lr \Hom(\one, \onehat(1)[2]) = \RR,$$ 
involving the motivic homology and the \emph{Arakelov motivic cohomology} of $M$ is a perfect pairing of finite-dimensional $\RR$-vector spaces.
This conjectural perfectness is very interesting in its own right. For example, special cases of it are equivalent to the Beilinson-Soul\'e vanishing conjecture (\refth{Beilinson.Soule}) and the Beilinson-Parshin conjecture (\refth{reductionperfFp}). It also allows to equivalently reformulate \refeq{pole.order} using the Euler characteristic of Arakelov motivic cohomology:
$$\ord_{s=0} L(M, s) = - \widehat \chi (M).$$
Most importantly, though, it allows to express the following conjecture for the special $L$-value $L^*(M, 0)$ up to a nonzero rational factor, using the determinants of the pairings $\pi_{M[i]}$: 
$$L^*(M, 0) \equiv \Pi_M^{-1} \mod {\Q^\x},$$
where 
$$\Pi_M := \prod_{i \in \Z} \det (\pi_{M[i]})^{(-1)^{i}} (\in \RR^\x / \Q^\x).$$
\xconj

The Arakelov motivic cohomology referred to above is a new cohomology established in 
\cite{HolmstromScholbach:Arakelov, Scholbach:ArakelovII} (or see \refsect{Arakelov}). 
It can be thought of as a cohomology with compact support, where ``compact'' refers to the compactification of $\SpecZ$. More precisely, it is characterized by a long exact sequence
$$\dots \r \Hhat^n(M) \r \H^n(M) \stackrel \ch \r \HD^n(M) \r \Hhat^{n+1}(M) \r \dots$$
involving the Chern class map $\ch$ (also known as the Beilinson regulator) between motivic cohomology and Deligne cohomology. 

This conjecture is related to existing conjectures on $L$-functions as follows:

\theo (see Theorems \ref{theo_summary.poles}, \ref{theo_summary} for the precise statements) 
Assuming the existence of the category of mixed motives (see \refax{mixed}), \refco{summary} is essentially equivalent to the conjunction of the conjectures \ref{conj_Beilinson}, \ref{conj_SouleII}, \ref{conj_Tate} of Beilinson, Soul\'e and Tate on special $L$-values of motives over $\Q$ and $\zeta$-functions \`a la Hasse-Weil of schemes over $\Z$ and over $\Fp$, respectively.
\xtheo

Recall that the subcategory $\DATM(\Z)$ of Artin-Tate motives is the triangulated subcategory generated by direct summands of motives of number rings $\OF$ and finite fields $\Fq$.
Only allowing $\Q$ and $\Fp$ instead of arbitrary $\OF$ and $\Fq$, we get the triangulated category $\DTM(\Z)$ of Tate motives. Note that these motives have rational coefficients.
These categories do enjoy a motivic t-structure whose hearts are denoted $\MATM(\Z)$ and $\MTM(\Z)$, respectively \cite{Scholbach:Mixed}. We get the following unconditional result:

\coro 
\mylabel{coro_ATM}
The perfectness of the Arakelov intersection pairing, as well as the pole order formula \refeq{pole.order} holds for any Artin-Tate motive over $\Z$.
The formula for the special $L$-value holds for all motives in the triangulated category generated by motives $\M(\OF)$ and $\M(\Fq)$, in particular for any Tate motive, i.e., any motive in $\DTM(\Z)$.
More generally, for any $M \in \DATM(\Z)$,
$$L^*(M, 0) \cdot \Pi_M$$
is a torsion element of $\RR^\x / \Q^\x$.

In particular, Beilinson's conjecture holds for any smooth projective variety $X_\eta / \Q$ such that $\h^j(X_\eta)$ is a mixed Tate motive ($j \in \Z$).
Examples of such varieties include linear varieties \cite[Section 14]{Jannsen:Mixed}, \cite{Totaro:Chow}, such as toric varieties and Grassmannians.
Similarly, Beilinson's conjecture holds up to the $m$-th root of a nonzero rational number if $\h^j(X_\eta)$ is a mixed Artin-Tate motive.
\xcoro

\pf
We first show that for any $M \in \DATM(\Z)$, there is some $m > 0$ such that $m M := M^{\oplus m}$ lies in the triangulated subcategory $L \subset \DATM(\Z)$ generated by motives of the form $\M(\OF)(n)[1]$ and direct factors of $\M(\Fq)$, for any $q=p^r$, $n \in \Z$ and any number ring $\OF$.
This statement is unrelated to the Arakelov intersection pairing and $L$-functions. 
It is enough to show this for $M$ being a direct summand of $\M(\OF)(n)[1]$.
By definition of $\eta_{!*}$, see \cite[Section 5.4]{Scholbach:fcoho}, $M' := \eta_{!*} \eta^* M$ lies in the triangulated category generated by $M$ and motives of the form $i_* N$, where $N \in \DATM(\Fp)$ and $i: \SpecFp \r \SpecZ$. As $i_* N \in L$ for all $N \in \DATM(\Fp)$, it is enough to show $m M' \in L$. 
Note that $M_\eta := \eta^* M[-1]$ is a direct summand of $\M(F)(n)$. After twisting by $-n$, these two motives are Artin motives over $\Q$ (with rational coefficients). This category is equivalent to continuous rational $\Gal(\Q)$-representations. For some finite quotient $G = \Gal (E / \Q)$ of $\Gal(\Q)$, $\M(F)$ and $M_\eta$ factor over $G$.
By Artin induction \cite[II.13.1, Th\'eor\`eme 30]{Serre:Representations}, there is an equality in $K_0(\Q[G])$, the $K_0$-group of the group ring of $G$ (with rational coefficients) $m [M_\eta(-n)] = \sum_i l_i [\ind_H^G \Q]$, where $m, l_i \in \Z$, $m > 0$, and $H$ runs over the cyclic subgroups of $G$.
The functor $\eta_{!*}[1]$ does not in general send a short exact sequence
$$E_\eta: \ 0 \r M_{\eta, 1} \r M_{\eta, 2} \r M_{\eta, 3} \r 0$$
in $\MATM(\Q)$ to a distinguished triangle in $\DATM(\Q)$. 
However, for a sufficiently small open $j: U \subset \SpecZ$, there is a similar short exact sequence $E_U$ in $\MATM(U)$ such that $\eta^* E_U [-1] = E_\eta$ and such that $\eta_{!*} M_{\eta, n}[1] = j_{!*} M_{U, n}$ for all $n$. As $j_!$ is triangulated, $j_! (E_U)$ is a distinguished triangle in $\DATM(\Z)$. 
Moreover, $j_{!*} M_{U, n}$ lies in a distinguished triangle whose other vertices are $j_! M_{U, n}$ and $i_* N$, where $i: Z \r \SpecZ$ is the complement of $j$ and $N \in \DATM(Z)$. Therefore, 
if $\eta_{!*} M_{\eta, j}[1] \in L$ for two out of the three $M_{\eta, j}$'s, it is true for the third. Noting that $\ind_H^G \Q$ corresponds to the motive $\M(E^H)$ of the subfield $E^H \subset E$ fixed by $H$ and $\eta_{!*} \M(E^H)[1] = \M(\mathcal O_{E^H})[1] \in L$, we obtain $m \eta_{!*} M_\eta[1] \in L$. 

For any number field $F$ and number ring $\OF$, the conjectured pole order formula, the special value and the perfectness of the Arakelov intersection pairings for $\M(\OF)(n)[1]$ are (unconditionally, by \refpr{comparison}, \refre{vanish.ATM}, and \refth{comparisonBeilinson}) equivalent to Beilinson's conjecture for $\M(F)(n) \in \MATM(\Q)$ which does hold by Borel's work \cite{Borel:Cohomologie}. 
The three conjectures also hold for direct factors of $\M(\Fq)$ by Quillen's computation of $K$-theory of finite fields \cite{Quillen:FiniteField}. By \refth{structure}, the three conjectures therefore hold for any motive in $L \subset \DATM(\Z)$.

Now, let $M \in \DATM(\Z)$ be any Artin-Tate motive. There is an $m > 0$ such that $mM \in L$.
Since the Arakelov intersection pairings are induced by the composition of morphisms in $\DMBei(\Z)$, the map $r_{mM} : \H_{-2}(mM, -1) \r \Hhat^0(mM)\dual$ induced by $\pi_{mM}$ is clearly the $m$-fold direct sum of the map $r_M$ induced by $\pi_M$. Hence the perfectness of $\pi_{mM}$, i.e., $r_{mM}$ being an isomorphism, implies the perfectness of $\pi_M$.
Moreover, we have $(L^*(M, 0) \Pi_M)^m  = L^*(mM, 0) \Pi_{m M} = 1 \in \RR^\x / \Q^\x$, i.e., $L^*(M, 0) \Pi_M$ is torsion in $\RR^\x / \Q^\x$. 
Similarly, $m (\ord_{s=0} L(M, s) + \chi(M\dual(-1))) = \ord_{s=0} L(mM, s) + \chi(mM\dual(-1)) = 0 \in \Z$, so that $\ord_{s=0} L(M, s) + \chi(M\dual(-1)) = 0$, i.e., the pole order formula holds.

The last statement follows immediately.
\xpf

\refco{summary} is compatible with the functional equation of $L$-functions. It is also stable under distinguished triangles (\refth{structure}). While the latter is a formal consequence of the setup, it is a key difference between our conjecture and Beilinson's conjecture for mixed motives over $\Q$. It allows to break up a motive into smaller pieces by means of distinguished triangles. This technique is unapplicable when working with Beilinson's original conjecture for motives over $\Q$. Moreover, \refco{summary} gives more freedom because it allows to work in the larger category of all geometric motives, as opposed to just smooth and projective varieties. 
It should be noted, though, that the proof of the equivalence of Beilinson's $L$-value formula and \refco{L} is formal, so that proving Beilinson's conjecture for any example not covered by techniques such as the ones in \refcor{ATM} will require new ideas.


The idea of reinterpreting the data in Beilinson's conjecture in terms of motives over $\Z$ is due to Huber. More precisely, a mixed motive $M_\eta$ over $\Q$ corresponds to the mixed motive $\eta_{!*} M_\eta[1]$ over $\Z$. This is reified for $L$-functions by \refeq{LQZ} and on the motivic side by an appropriate interpretation of $f$-cohomology \cite{Scholbach:fcoho}. The non-multiplicativity of $L$-functions (cf.\ \refbe{propertiesL}) is explained by the failure of $\eta_{!*}$ to be exact. $L$-functions of motives over $\Z$ \emph{are} multiplicative, though. 

This non-multiplicativity, which is a heavy technical burden, has been addressed by Scholl by introducing a category $\MM(\Q / \Z)$ of mixed motives over $\Z$ \cite{Scholl:Remarks} (different from the one used here) by imposing non-ramification conditions. The (conjectural) value of the groups $\Ext^a_{\MM(\Q / \Z)} (\one, \h^{b-1}(X_\eta, m))$ is closely related to the computation of $\H^{*}(\eta_{!*} \h^{-b+1}(X_\eta, -m) [1])$ (\refth{summaryH1f}). As for the special $L$-values, a conjecture of Scholl \cite[Conj. C]{Scholl:Remarks} says that some $M_\eta \in \MM(\Q / \Z)$ is critical (i.e., its period map is an isomorphism, equivalently all weak Hodge cohomology groups $\Hw^*(M_\eta)$ vanish) if 
$$\Ext^a_{\MM(\Q / \Z)} (M_\eta, \one(1)) = \Ext^a_{\MM(\Q / \Z)} (\one, M_\eta) = 0 \text{ for }a = 0, 1.$$
Moreover, a reduction technique transforming any motive $M_\eta$ into one satisfying these vanishings is given, so that Deligne's conjecture \cite[Conj. 2.8.]{Deligne:Valeurs} concerning the $L$-value of critical motives can be applied. In similar spirit, the non-multiplicativity of $L$-functions of motives over $\Q$ has been addressed by Fontaine and Perrin-Riou by introducing the notion of $f$-exact sequences, which are ones where one does save multiplicativity \cite[III.3.1.4]{FPR}. However, such exact sequences seem to be hard to characterize. The formulation of \refco{summary} resembles their approach; for example the pole order in \ocs is expressed as an Euler characteristic of $f$-cohomology. Using a ``cohomology with compact support'' to predict special $L$-values was already suggested by Beilinson \cite[5.10.F]{Beilinson:Height}. The category of motives over $\Z$ is both the appropriate home for this idea and allows for the strikingly compact and beautiful formulation of the $L$-values conjecture by 
overcoming the technical obstacles related to motives over $\Q$.

The idea to recast special $L$-values of motives as determinants of appropriate pairings was explored by Deninger and Nart \cite{DeningerNart}, who show that the motivic height pairing of \cite{Scholl:Height} can be represented by concatenating morphisms in the derived category of an appropriate category of motives. 

\refco{summary} is the first conjecture that predicts the special values of $\zeta(X)$ modulo $\Q^\x$ at all places ($X / \Z$ regular projective; see \refbs{specialHasseWeil}). A reformulation of the Tamagawa number conjecture in terms of the Weil-\'etale cohomology due to Flach and Morin predicts the special value of $\zeta(X)$ at $s=0$ up to sign \cite[Prop. 9.2]{FlachMorin:Weil}. It remains to explicitly compare the compatibility of the approach of \ocs and \refco{summary}. I expect that similar techniques as the ones in this paper allow to refine \refco{summary} to a conjectural $L$-values formula, up to sign, at all places. However, this remains to be done. 
  
This paper has its origins in a part of my PhD thesis. It is a pleasure to thank Annette Huber for her advice during this time. I thank Andreas Holmstrom for the collaboration on Arakelov motivic cohomology \cite{HolmstromScholbach:Arakelov}. I also thank Denis-Charles Cisinski, Fr\'ed\'eric D\'eglise and Bruno Kahn for helpful conversations. 
Finally, I thank the referee for suggesting many improvements to this article, in particular concerning the formulation of \refcor{ATM}.

\section{Preliminaries}

\subsection{Determinants and $\Q$-structures} \mylabel{sect_LA}

For any ring $R$, let $\underline R$\index{R@$\underline{R}$} be the category of finitely generated $R$-modules. Let $K$ be a field. The \Def{determinant} $\det V$\index{det@$\det$} of $V \in \underline K$ is $\det V := \Lambda^{\dim V} V$. Its $K$-dual is denoted $\det^{-1} V$. For $V_* \in \D^\bound(\underline K)$, the derived category, we set $\det V_* := \bigotimes_i \det^{(-1)^i} \H^i(V_*).$ We abbreviate $\det H^* := \det^{(-1)^i} H^i$ for some $H^i \in \underline K$, $i \in \Z$.

Let $A, B \in \underline \Q$ and let $f: A_\RR \r B_\RR$ be an $\RR$-linear map. We do not assume that it respects the rational subspaces. The ``usual'' determinant of $f$, which is well-defined up to a nonzero rational factor agrees, modulo $\Q^\x$ with the image of $1$ under the map $\Q \cong \det A \t \det^{-1} B \r \det A_\RR \t \det^{-1} B_\RR \cong \RR$. Here the right hand isomorphism is induced by $f$.

A complex with \emph{$\Q$-structure}\index{Q-structure@$\Q$-structure} is a complex $V_*$ of $\RR$-vector spaces that is quasi-isomorphic to one in $\D^\bound(\underline \RR)$ together with a non-zero map of $\Q$-vector spaces $d_{V_*}: \Q \r \det V_*$. In concrete situations, we usually have a distinguished identification $\det V_* \cong \RR$. In that case, we may also call $\det V_*$ the real number that is the image of $1 \in \Q$ under $d_{V_*}$ and the given identification.

Maps of complexes with $\Q$-structures are usual maps of complexes; they are \emph{not} required to be compatible with the map $d_{V_*}$. For a map $f: V_* \r W_*$ of complexes with $\Q$-structures the cone of $f$ is endowed with the following $\Q$-structure:
$$\Q \stackrel{d_W \t (d_V)^{-1}} \lr \det W_* \t \det^{-1} V_* \cong \det \cone(f).$$

Define a category $\DbRQdet$\index{DbRQdet@$\DbRQdet$} to consist of such complexes. Its morphisms are given by maps of complexes up to quasi-isomorphism (not necessarily respecting the $\Q$-structures). We say that a triangle $A \r B \r C$ of objects in $\DbRQdet$ is \Def{multiplicative}\index{multiplicative triangle} if it is distinguished in $\D^\bound(\underline \RR)$ after forgetting the $\Q$-structure and $\det B = \det A \det C$ in the sense that the following diagram (whose right hand isomorphism stems from the triangle) is commutative:
$$\xymatrix{
\Q \ar[r]^{d_C} \ar[dr]_{(d_A)^{-1} \t d_B} & \det C \ar[d]^\cong \\
& \det^{-1} A \t \det B.
}$$



\subsection{Motives} \mylabel{sect_motives}

Our work takes place in the category $\DMBei(S)$ of \emph{Beilinson motives} over $S$, where $S$ is either a finite field, a number ring $\OF$, or a number field $F$. Cisinski and D\'eglise defined this category to be an appropriate subcategory of Morel and Voevodsky's stable homotopy category $\SH(S)_\Q$ (with rational coefficients) \cite{CisinskiDeglise:Triangulated}. The category $\DM_{\Beilinson}(S)$ is tensor-triangulated, $\Q$-linear, and closed under arbitrary direct sums. Its tensor unit is denoted $\one_S$ or just $\one$. Given some scheme $f: X \r S$ (always tacitly supposed to be separated and of finite type), the \Def{motive} of $X$ and the \Def{motive with compact support} are defined as 
\eqn \mylabel{eqn_motiveX} 
\M(X) := f_! f^! \one_S, \ \ \Mc(X) := f_* f^! \one_S.
\xeqn
Here $f_!: \DMBei(X) \r \DMBei(S)$ etc. are the functors defined in \ocs This determines a covariant functor $\M: \Sch / S \r \DMBei(S)$ and likewise, but just for proper maps, with $\Mc$. The motive of the projective line decomposes as $\M(\P) = \one \oplus \one(1)[2]$. In $\DMBei(S)$, tensoring with $\one(1)[2]$ is invertible and we write $M\twi n := M \t (\one (1)[2])^{\t n}$ for any $n \in \Z$. 
For $M \in \DMBei(S)$, we put
\eqnarr
\mylabel{eqn_motivichomology} 
\H_i(M, p)& := & \Hom_{\DMBei(S)}(\one(p)[i], M), \\
\H^i(M, p) &:= & \Hom_{\DMBei(S)}(M, \one(p)[i]).
\xeqnarr
For a regular base $S$ and a regular, projective or affine (but not necessarily flat) scheme $X$ over $S$ and $G := \M(X)(-m)$, \emph{motivic cohomology} of $X$ is given by
$$\H^i(X, m) := \H^i(G) = \Hom_S(f_! f^! \one(-m), \one[i]) = \Hom_X(\one, \one(m)[i]) = K_{2m-i}(X)_\Q^{(m)},$$
using the purity isomorphism $f^! \one_S = f^* \one_S \twi d = \one_X \twi d$, where $d = \dim X - \dim S$. As a consequence of resolution of singularities, the full subcategory $\DMBeic(S) \subset \DMBei(S)$ of compact objects agrees with the thick subcategory generated by such motives $G$, for any base $S$ as above.\mylabel{compactobjects} We refer to objects of $\DMBeic(S)$ as \emph{geometric motives over $S$}. For a perfect field $S$, there is a natural equivalence of categories \cite[Theorem 15.1.4]{CisinskiDeglise:Triangulated}
\eqn \mylabel{eqn_DMgm}
\DMBeic(S) \stackrel \cong \lr \DMgm(S)_\Q
\xeqn  
with Voevodsky's triangulated category of geometric motives (with rational coefficients) \cite{Voevodsky:TCM}. It sends the motive $\M(X) \in \DM_{\Beilinson, c}(\Q)$ of a smooth $S$-scheme in the sense of \refeq{motiveX} to the motive $\Mgm(X)$ of $X$ in Voevodsky's sense.

The category $\DMBeic(S)$ is equipped with a notion of \Def{weight}: there are full (non-triangulated) subcategories $\DMBeic^{\wt \leq n}(S)$, $\DMBeic^{\wt \geq n}(S)$ such that $f_! \one (a)[2a+n]$ lies in the subcategory 
\eqn
\mylabel{eqn_weights}
\DMBeic^{\wt = n}(S) := \DMBeic^{\wt \leq n}(S) \cap \DMBeic^{\wt \geq n}(S)
\xeqn 
of objects of pure weight $n$, for all $a, n \in \Z$ and all proper maps $f: X \r S$ with regular domain $X$ \cite{Bondarko:Weight, Hebert:Structures}. For any map $f$ (of finite type), the functors $f_!$, $f^*$ preserve the subcategories $\DMBeic^{\wt \leq n}(-)$ and dually for $f^!$, $f^*$. 

The dual of any geometric motive $M$ is defined as $M\dual = \IHom_{\DMBei(S)}(M, \one)$. Dualizing exchanges $!$ and $*$: for example, for any map $f$, $(f_! f^! \one)\dual = f_* f^* (\one\dual)$ which is canonically isomorphic $f_* f^* \one$. Therefore, the natural map 
\eqn \mylabel{eqn_reflexive}
M \r (M\dual)\dual
\xeqn
is an isomorphism for any $M \in \DMBeic(\Z)$ \cite[14.3.31]{CisinskiDeglise:Triangulated}. This yields a canonical isomorphism $\H_0(M\dual, 0) = \H^0 (M, 0)$.\mylabel{VerdierDual}

\defi
\mylabel{defi_generically.smooth}
Let $S \subset \SpecOF$ be an open subscheme. A motive $M \in \DMBeic(S)$ is called \Def{smooth} if the natural map \cite[Section 2.3.2]{Ayoub:Six1} \footnote{The use of this canonical map, as opposed to a mere noncanonical isomorphism, was suggested by Bruno Kahn.} 
$$i^*M \twi {-1} = i^*M \t i^! \one \r i^! (M \t \one) = i^! M$$
is an isomorphism for all closed points $i: \SpecFpp \r S$. 
A motive $M \in \DMBeic(S)$ is \Def{generically smooth} if $j^* M$ is smooth for some open subscheme $j: U \subset S$. 
\xdefi

Since $\M(X)(m)$ is smooth provided $X/S$ is smooth and proper, every motive $M \in \DMBeic(S)$ is generically smooth.
We write $\eta: \Spec F \r \SpecOF$ for the generic point. 

In order to interpret Beilinson's conjecture for mixed motives over $\Q$ in terms of motives over $\Z$ we need to assume the conjectural framework of mixed motives over $F$, $\SpecOF$ and $\Fq$. The precise axioms we are staking on are listed in \cite[\fcohoAxiomSections]{Scholbach:fcoho}, so we only summarize them briefly and refer to \lcs for a more complete discussion.
\footnote{Unlike this paper, \ocs is written with a contravariant notation of motives. This induces a number of changes in notation: every $f_!$, $f^!$ gets replaced by a $f_*$ and $f^*$, and vice versa. Moreover, a twist and shift $M(m)[n]$ corresponds to $M(-m)[-n]$ here. Both here and there, the normalization of the $t$-structure is such that $\one[1] \in \MM(\OF)$, while $\one \in \MM(F)$.}
\footnote{
The decomposition axiom for smooth projective varieties formulated in \cite[Axiom 4.13]{Scholbach:fcoho} is not needed: it is only used in \cite[Lemma 5.10]{Scholbach:fcoho} to show that a certain motive is generically smooth, but this is inconditionally true for any motive by the remark after \refde{generically.smooth}.} 
Note that the corresponding statements for the triangulated category of \Def{Artin-Tate motives} $\DATM(\OF) \subset \DMBeic(\OF)$, which is the triangulated subcategory of $\DMBei(\OF)$ generated by $\M(V)(m)$ where $V \r \SpecOF$ is a quasi-finite, but not necessarily flat map and $m \in \Z$, have been shown in \cite{Scholbach:Mixed}. 

\axio
\mylabel{axio_mixed}
\begin{enumerate}[(i)]
\item 
\cite[Axioms 4.1, 4.2]{Scholbach:fcoho}
$\DMBeic(S)$ is conjectured to enjoy a non-degenerate $t$-structure whose heart $\MM(S)$ is called the category of \Def{mixed motives}. The cohomological dimension of $\MM(S)$ is conjectured to be $0$ ($S = \Fq$) and $1$ ($S = F$), respectively. The truncation with respect to the $t$-structure is 
denoted $\pH^*$. We write $\h^i(X, n)$ for $\pH^i(\M(X)(n))$. The $t$-structures are normalized by declaring $\one \in \MM(S)$ when $S = F$, $\Fq$ and $\one[1] \in \MM(\OF)$, respectively. For example, $\h^{-1}(\P_{\OF}) = \one_{\OF}(1)[1]$, $\h^{-2}(\P_F) = \one_F(1)$. More generally, $\eta^* [-1]$ is $t$-exact and $\eta^* \h^{-b}(X, -m) = \h^{-b-1}(X_\eta, -m)$ for any scheme $X / \OF$ with generic fiber $X_\eta$. 
\item 
\mylabel{item_exactness}
\refaxtstructureRealizations \ 
The key requirement on the $t$-structure is that realization functors of the form $\DMBeic(S) \r \D^\bound(\C)$ are to be exact (see \lcs \ and around \refeq{exactnessell} for the $\ell$-adic realization over $\Z[1/\ell]$). 
In the guise of a spectrum representing the cohomology theory, the exactness requirement is to be understood as in \refeq{exactnessBetti}.
\item
\cite[Axioms 4.4, 4.6, 4.11]{Scholbach:fcoho} 
Any mixed motive is conjectured to have a weight filtration which is compatible with the weight formalism mentioned around \refeq{weights}. The pure objects in $\MM(K)$ (for any field $K$) are conjectured to be identified with the category $\PureMot_\num$ 
of pure motives with respect to numerical equivalence. This implies that the pure objects in $\MM(K)$ form an abelian semi-simple category \cite[Th. 1]{Jannsen:Motives}. Moreover, homological and numerical equivalence are conjectured to agree. The cohomology functors $\pH^*$ belonging to the motivic $t$-structure are supposed to respect the weights, i.e., given some $M \in \DMBeic^{\wt=w}$, $\pH^n (M) 
\in \MM$ is pure of weight $w+n$. For example, for a smooth projective scheme $X / S$, $\M(X)(-m) = f_! f^! \one(-m) \in \DMBeic^{\wt=2m}$, so that $\h^{-b}(X, -m)$ is pure of weight $2m-b$. Morphisms of mixed motives are expected to respect weights strictly, thereby giving constraints on the existence of maps between motives.
\end{enumerate}
\xaxio

\emph{\noindent In the remainder of this paper we assume that the axioms concerning mixed motives over open subschemes of $\SpecOF$, $\Fq$ and $F$ hold.} \\

Given a mixed motive over $\Q$, $M_\eta \in \MM(\Q)$, pick any $M \in \MM(\Z)$ satisfying $M_\eta = \eta^* M [-1]$ and some open subscheme $j: U \r \SpecZ$ such that $j^* M$ is smooth. We call
$$\eta_{!*} (M_\eta [1]) := j_{!*} j^* M := \im (j_! j^* M \r j_* j^* M) \in \MM(\Z)$$
the \Def{generic intermediate extension} of $M_\eta[1]$. This is explained and shown to be well-defined in \cite[\fcohoSectionGenericIntermediateExtension]{Scholbach:fcoho}. 
We apply this to $M_\eta = \h^{-b-1}(X_\eta, -m)$ and $M = \h^{-b}(X, -m)$, where $X_\eta / \Q$ is smooth projective and $X / \Z$ is any projective (not necessarily regular) model of $X_\eta$ of constant dimension $d$. Throughout this paper, we write 
\eqn \mylabel{eqn_E}
E := \eta_{!*} \eta^* \h^{-b}(X, -m) = \eta_{!*} (\h^{-b-1}(X_\eta, -m)[1]) \in \MM(\Z).
\xeqn 
This motive is pure of weight $w := 2m-b$. Its motivic cohomology is given by the following theorem:

\theo \mylabel{theo_summaryH1f}
With the above notation, we write $\H^b(X_\eta, m)_\Z := \im (\H^b(X, m) \r \H^b(X_\eta, m))$. Moreover, let $\CH^m(X_\eta)_{\Q,\hom}$ be the subgroup of the Chow group of cycles homologically equivalent to zero and $\CH^m(X_\eta)_\Q / \hom$ the group of cycles modulo homological equivalence (tensored with $\Q$). Then
$$\H^a(E) = \H^a(\eta_{!*} \h^{-b-1}(X_\eta, -m) [1]) = \left \{ \begin{array}{ll}
\CH^{m}(X_\eta)_\Q / \hom  & a = 1,\, w=1 \\
0 & a=1,\, w \neq 1 \\
0 & a=2,\, w \leq 1 \\
\CH^{m}(X_\eta)_{\Q,\hom} & a = 2,\, w= 2 \\
\H^{b+2} (X_\eta, m)_\Z & a = 2,\, w \geq 3 \\
0 & a = 3,\, w \leq 2 \\
? & a = 3, \, w \geq 3 \\
0 & a > 3, a < 1
\end{array}
\right.
$$
\xtheo

\pf
Everything except the cases $a=2$, $w\leq 1$ and $a=3$, $w \leq 2$ is shown in \cite[\fcohoLemmcohomdimOF, \fcohoTheoSummaryHOneF]{Scholbach:fcoho}. For $a=2$ and $w \leq 1$, the map 
$$\H^2(E) \r \H^2(\eta^* E) = \H^{1}(\h^{-b-1}(X_\eta, -m)) \r \H^{b+2}(X_\eta, m) = \CH^m(X_\eta, w-2) = 0$$ 
is injective: for the first map this is \cite[\fcohoLemmprepinjectivity]{Scholbach:fcoho}, the second one is because the cohomological dimension of $\DMBeic(\Q)$ is one \refaxcohomdim. For $a = 3$, $w \leq 2$, we use the exact localization sequence
$$\dots \r \oplus_{p} \H^3(i_p{}_* i_p^* E) \r \H^3(E) \r \H^3(\eta^* E) = \H^2(\eta^*[-1] E) = 0.$$
The right hand vanishing is again because the cohomological dimension of motives over $\Q$ being one. Also by cohomological dimension we have
$$\H^3(i_p{}_* i_p^* E) = \Hom (i_p^* E, i_p^! \one[3]) = \Hom (i_p^* E, \one(-1)[1]) = \Hom_{\MM(\Fp)}(\pH^{-1} i_p^* E(1), \one).$$ 
The functor $i^*$ preserves negative weights, i.e., $\wt (\pH^{-1} (i_p^* E(1))) \leq \wt (E) - 1 - 2 = w-3$. By strictness of the weight filtration the group therefore vanishes for $w \leq 2$.
\xpf

In accordance with \refco{motcomp} (see the case $w\leq1$ in the proof of \refpr{comparison}) I expect $\H^3(E) = 0$ for arbitrary weight $w$. See the introduction for the relation of this to Scholl's notion of mixed motives over $\Z$.
%
For Artin-Tate motives, the expected vanishing holds unconditionally for all weights:   
%

\theo \mylabel{theo_H1Tate}
Let $M_\eta$ be an Artin-Tate motive over $F$, concentrated in cohomological degree $-1$. Then $\H^3(\OF, \eta_{!*} M_\eta) = 0$.
\xtheo

\pf
There is some $j: U \subset \SpecOF$ and a smooth Artin-Tate motive $M \in \MATM(U) = \MM(U) \cap \DATM(U)$ such that $M_\eta = \eta^* [-1] M$. Shrinking $U$ further (using $j'_{!*} j'^* M \cong M$ for some $j': U' \subset U$, as $M$ is smooth), we may assume by the standard splitting routine \cite[\mixedLemmSplit]{Scholbach:Mixed} that there is an etale Galois cover $f': V' \r U$ such that $f'^* M$ is a mixed Tate motive over $V'$. The map $M \r f'_* f'^* M \stackrel{\cong}\leftarrow f'_! f'^! M \r M$ is $\deg f' \cdot \id_M$, so $M$ is a direct summand of $f'_* f'^* M$, since we use rational coefficients.  
The functor $f'_* = f'_!$ preserves Artin-Tate motives and is exact \cite[\mixedTheoExactness]{Scholbach:Mixed}.  Hence $j_{!*} f'_* f'^* M = f_* j'_{!*} f'^* M$. Here $f: V \r \SpecOF$ is the normalization of $\OF$ in the function field of $V'$ and $j': V' \r V$ is the corresponding open immersion. Consequently, 
$$\H^3(\OF, \eta_{!*} M_\eta) = \H^3(\OF, j_{!*} M) \subset \H^3 (V, j'_{!*} f'^* M) = \Hom_V(j'_{!*} f'^* M, (\one[1])[2]) = 0,$$ 
since the cohomological dimension of mixed Tate motives over $V$ is one, as opposed to two for Artin-Tate motives \cite[\mixedSatzCohomdimAT]{Scholbach:Mixed}.
\xpf


The following conjecture will be needed to deal with motives over $\Fp$. 

\conj (Beilinson) \mylabel{conj_numrat} 
Let $X / \Fq$ be smooth and projective. Up to torsion, numerical and rational equivalence agree on $X$.
\xconj

Recall that homological equivalence lies between these two equivalence relations \cite[3.2.1]{Andre:Motifs}, so under \ref{conj_numrat}, all three agree. The second important consequence of \ref{conj_numrat} is that the category of pure Chow motives over $\Fq$ is semisimple by Jannsen's theorem.

To study $L$-functions, we need some $\ell$-adic realization functor. We use the machinery developed recently by Ayoub \cite{Ayoub:Realisation}. It allows the base scheme to be $\Z[1/\ell]$. Let $\ell$ be an odd prime number and $S$ a scheme that is of finite type over $\Z$ or $\Q$ such that $\ell$ is invertible on $S$ (cf.\ \cite[Hyp. 6.5]{Ayoub:Realisation}). Define the $\ell$-adic realization functor as the following composition 
\eqnarra
\mylabel{eqn_ell-adic}
(-)_\ell : \DM_{\Beilinson}(S) \stackrel{F_1} \r \SH(S)_\Q & \stackrel{F_2} \r & DA^\et(S, \Ql) \\
& \stackrel{R_\ell}\r &\D(\Shvv_\et(S, \Ql)) \stackrel{F_3}\lr \D(\Shvv_\et(S, \Ql)) \nonumber
\xeqnarra
The functor $F_1$ is the inclusion of the category of $\one$-modules in $\SH(S)_\Q$. The category $DA^\et(S, \Ql)$ is the homotopy category of the model category of symmetric $\P$-spectra of complexes of $\ell$-adic presheaves on $\Sm/S$, endowed with the $\A$-\'etale-local model structure. The functor $F_2$ is obtained by combining the natural free abelian group functor $\Delta^\op \Sets \r \Com(\Ab)$ and the sheafification (from Nisnevich sheaves to etale sheaves), see e.g. \cite[5.3.28, 5.3.37]{CisinskiDeglise:Triangulated}. The functor $R_\ell$ is Ayoub's $\ell$-adic realization functor. We append the contravariant functor $F_3: M \mapsto \IHom(M, f^! \Ql)$, where $f: S \r \Spec \Z$ is the structural map (and $\IHom$ denotes the derived inner $\Hom$). For any map $g: X \r Y$ of quasi-projective $S$-schemes, the functors $F_1$, $F_2$, $R_\ell$ commute with $g_!$, $g_*$, $g^*$ and $g^!$ and, when applied to compact objects, with $\IHom$ \cite[Thm. 6.6]{Ayoub:Realisation}. Finally, $F_3$ exchanges $!$ and $*
$, e.g. $F_3 (g^* M) = g^! F_3(M)$ for $M \in \D(\Shvv_\et(S, \Ql))$. Therefore, for some quasi-projective scheme $f: X \r S$, $(\M(X)(-m)[-n])_\ell = f_* f^* \Ql(m)[n]$. This property is also satisfied for Huber's and Ivorra's realization functors provided $S$ is a field \cite{Huber:Realization, Ivorra:Realisation}. Thus, for the mere definition in \ref{defi_LfunctionoverOF}, these realization functors are sufficient, but \refle{excreal} relies on a realization functor over $\Z[1/\ell]$.

\mylabel{exactnessell} The exactness requirement for the functor $-_\ell$ mentioned in \refax{mixed}\refit{exactness} means that the restriction of $-_\ell$ to $\DMBeic(S)$ is exact with respect to the (conjectural) motivic $t$-structure and the $t$-structure on $\D(\Shvv_\et(S, \Ql))$ (which is the obvious one if $S$ is a field and the perverse $t$-structure for $S = \SpecZ[1/\ell]$, see \cite[\fcohoSectionThree]{Scholbach:fcoho}. For example, for a quasi-projective variety $X$ over a field it implies 
\eqn
\mylabel{eqn_exactnessell}
(\h^{-b}(X))_\ell = \H^b(X, \Ql).
\xeqn 

\section{Arakelov motivic cohomology}

\subsection{Deligne cohomology} \mylabel{sect_absolute}

A key input to Beilinson's conjecture \ref{conj_Beilinson} is Deligne cohomology. We recall its classical definition and the well-known interpretation in terms of weak Hodge cohomology. Then, we recall from \cite{HolmstromScholbach:Arakelov} the Deligne cohomology spectrum $\HD$ which is crucial for the definition of Arakelov motivic cohomology. In order to establish the $\Q$-structure on the groups represented by $\HD$, we explain how to apply the construction in \lcs to obtain spectra representing Betti and de Rham cohomology.  

Let $\an: \Sm/\CC \r \Sm^\an$ be the functor that associates to any smooth $\CC$-scheme the underlying complex analytic manifold. We also consider $\an: \Sm / \Q \text{ (or }\Sm / \RR\text) \r \Sm^{\an, G}$, where the target category consists of complex analytic manifolds with a $G$-action, $G := \Gal(\CC/ \RR)$. In this section, $X$ is a smooth scheme over $\Q$. We usually write $X^\an := \an(X)$ and $\Fr_\infty: X^\an \r X^\an$ for the conjugation. We also pick a smooth proper compactification $j: X \r \ol X$ (over $\Q$) such that $D := \ol X \backslash X$ is a divisor with strict normal crossings. We write $\Omega^*_{\ol X^\an}(\log D^\an)$ for the complex of meromorphic forms on $\ol X$ that are holomorphic on $X \subset \ol X$, and have at worst logarithmic poles at the divisor $D$. This complex is endowed with the Hodge filtration $F^p := \sigma_{\geq p}$, which is simply the brutal truncation. The variant using algebraic (i.e., K\"ahler) differential forms is denoted $\Omega^{*, \alg}_{\ol X}(\log D)$.
 The $C^\infty$-variant is denoted $E^*_{\ol X^\an}(\log D^\an)$. The subspace of real-valued forms is denoted $E^*_{\RR, \ol X^\an}(\log D^\an)$. These complexes are filtered by $F^p E^n_{\ol X^\an}(\log D^\an) = \oplus_{a+b=n, a \geq p} E^{a,b}_{\ol X^\an}(\log D^\an)$. To get rid of the choice of $\ol X$, put 
$$E^*(X) := \varinjlim_{\ol X}  E^*_{\ol X^\an}(\log D^\an),$$
and similarly for $E^*_\RR(X)$, $\Omega^*(X)$, $\Omega^{*, \alg}(X)$. Here, the colimit runs over the directed category of all compactifications $\ol X$ as above. Finally, let $\RR(p) := (2 \pi i)^p \RR \subset \CC$ be the constant sheaf. 

\defi \mylabel{defi_Delignecomplex}
Set $\RR_{\Deligne, D, \ol X}(p):= \cone (\R j_* \RR(p) \oplus F^p \Omega^*_{\ol X^\an}(\log D^\an) \r \R j_* \Omega^*_{X^\an})$. For example, if $X$ is proper, $\RR_\Deligne(p) \cong [\RR(p) \r \Omega^0_{X^\an} \r \dots \r \Omega^{p-1}_{X^\an}]$, with the terms lying in degrees $0$ to $p$. \Def{Deligne cohomology} of $X$ is defined as the $G$-invariant subspace of a sheaf hypercohomology group,
$$\HD^n(X, p) := \HH^n(\ol X^\an, \RR_{\Deligne, D, \ol X}(p))^G.$$
(The $G$-action is obtained by letting $G$ act on $\RR(p)$ as $a \mapsto \ol{\Fr_\infty^*(a)}$ and on $\Omega^*$ by $\omega \mapsto \Fr_\infty^*(\omega)$. This group does not depend on the choice of $\ol X$ \cite[Lemma 2.8]{EsnaultViehweg}.)
\xdefi

By definition, there is a long exact sequence
$$\dots \r (\HdR^i (X^\an) / F^p \H^i(X^\an, \Omega^*_X)^G \r \HD^{i+1}(X, m) \r \H^{i+1}(X^\an, \RR(m))^{(-1)^m} \r \dots.$$
Here the superscript denotes the $(-1)^m$-eigenspace of the $\Fr_\infty$-action on Betti cohomology of $X^\an$. This sequence induces an isomorphism 
\eqn \mylabel{eqn_Qstructure}
\det \HD^*(X, m) = \det^{-1} (\HdR^*(X^\an) / F^m)^G \t \det \H^* (X^\an, \RR(m))^{(-1)^m}.
\xeqn
The right hand side carries a natural $\Q$-structure stemming from the isomorphisms $\H^*(X^\an, \RR(m))=\H^*(X^\an, \Q(m)) \t_\Q \RR$ and $\H^*(\ol X^\an,  F^* \Omega^*_{\ol X^\an}(\log D^\an))^G \cong \H^*(\ol X_\RR, F^* \Omega^{*,\alg}_{\ol X_\RR}(\log D_\RR)) = \H^*(\ol X, F^* \Omega^{*,\alg}_{\ol X}(\log D)) \t_\Q \RR$ (GAGA). We use the above isomorphism to carry over the $\Q$-structure to the left hand side.
      
If $X$ is (smooth and) proper, the degeneration of the Hodge-de Rham spectral sequence and weight reasons give us short exact sequences (\lc)
\eqn
\mylabel{eqn_longDeligne}
0 \r \H^i (X^\an, \RR(m))^{(-1)^m} \r \H^i_\dR (X_\RR) / F^m \r \HD^{i+1}(X, m) \r 0
\xeqn
for $i - 2m \leq -2$ and, for $i-2m \geq 0$, 
\eqn 
\mylabel{eqn_longDeligne2}
0 \r \HD^i(X, m) \r \H^i (X^\an, \RR(m))^{(-1)^m} \r \H^i_\dR (X_\RR) / F^m \r 0,
\xeqn
respectively. In this case, each individual Deligne cohomology group carries a $\Q$-structure, as opposed to the general case of a merely smooth $X / \Q$.

Now, we recall Beilinson's notion of weak absolute Hodge cohomology. It is relevant to us because of its relation to archimedean factors of $L$-functions, see \refeq{ordLinfty}. It is based on Deligne's abelian category $\MHS_\Q(\RR)$ of mixed Hodge structures \cite[2.3.1]{Deligne:Hodge2}. The subscript $\Q$ indicates that we are considering $\Q$-vector spaces, ''$(\RR)$''  means that the structure is endowed with an action of $G=\Gal(\CC/\RR)$. For example, $\one(n)$ is the one-dimensional $\Q$-space, such that it is pure of weight $-2n$ and the Hodge filtration is concentrated in degree $-n$, and the non-trivial element of $G$ acts as multiplication by $(-1)^n$. Let 
$$\Com^\bound_\mathrm{H} = \{ C = (C_\dR, C_\Betti, C_c, i_\dR, i_\Betti) \}$$
be the category of bounded \Def{Hodge complexes} \cite[3.2]{Beilinson:Notes}. Its objects consist of a bounded bifiltered complex of $\Q$-vector spaces $(C_\dR, W_*, F^*)$, a filtered complex of $\Q[G]$-modules $(C_\Betti, W_*)$  and a filtered complex of $\CC$-modules with $\CC$-antilinear $G$-action, $(C_c, W_*)$, a filtered $G$-equivariant quasi-isomorphism $i_\Betti: (C_\Betti, W_*) \t_\Q \CC \r (C_c, W_*)$ ($G$ acts on the left hand term by the action on $C_\Betti$ and complex conjugation on $\CC$) and finally a filtered $G$-equivariant quasi-isomorphism $i_\dR: (C_\dR, W_*) \t_\Q \CC \r (C_c, W_*)$ (on the left, $G$ acts by conjugation on $\CC$). These data are subject to the requirement that the cohomology quintuple $\H^i(C)$ defined by the cohomologies of the various complexes and comparison maps has to be an object of $\MHS_\Q(\RR)$. Morphisms in the category $\Com^\bound_\mathrm{H}$ are required to respect the filtrations and the comparison quasi-isomorphisms. To any Hodge complex, we can associate 
its \Def{weak Hodge cohomology} \cite[3.13]{Beilinson:Notes} 
$$\RGw(C) := \cone[-1] \left (C_\Betti^G \t \RR \oplus F^0 C_\dR \t \RR \stackrel{i_\Betti - i_\dR}\lr C_c^G \right) \in \Com(\underline \RR).$$
This descends to a functor
$$\RGw: \DbH := \Com^\bound_\mathrm H / \text{quasi-isomorphisms} \r \DbRQdet.$$
Indeed, taking $G$-invariants and applying the Hodge filtration are exact operations, since morphisms of Hodge structures respect the Hodge filtration strictly \cite[2.3.5(iii)]{Deligne:Hodge2}. The $\Q$-structure on $\RGw(C)$ is the one stemming from the very definition, where $C_c^G$ is endowed with a $\Q$-structure using the one on $C_\dR$ via $i_\dR$. Set $\Hw^i(C) := \H^i(\RGw(C))$. 
A spectral sequence argument yields an exact sequence:
\eqn \mylabel{eqn_specweak}
0 \r \Hw^1 (\H^{i-1} C) \r \Hw^i (C) \r \Hw^0 (\H^i C) \r 0.
\xeqn
\ifpreprint
For $V \in \MHS_\Q(\RR)$, there is an exact sequence \cite[5.4]{Fontaine:Valeurs}
\small
$$0 \r \Hw^0(W_0 V) \r \Hw^0(V) \r \Hw^0(V / W_0 V) \r \Hw^1 (W_0 V) \r \Hw^1(V) \r \Hw^1(V / W_0 V) = 0.$$
\normalsize
\fi
Unlike absolute Hodge cohomology, i.e., the derived functor of $V \mapsto \Gamma_\MHS(V) := \Hom_\MHS(\one, V) = \Hw^0(W_0 V)$, the weak variant has a duality: the natural pairing (induced by $A \x A\dual \r \RR$ for any $\RR$-vector space $A$), 
\eqn \mylabel{eqn_dualityweak}
\Hw^i (C) \x \Hw^{1-i} (C\dual(1)) \r \Hw^1(\one(1)) = \RR,
\xeqn
is perfect for all $i$ \cite[Prop.III.1.2.3]{FPR}. 

The following well-known lemma states that weak Hodge cohomology is the same as Deligne cohomology. Recall the Hodge complex $\underline{\RG} (X, m)$ of \cite[Section 4]{Beilinson:Notes} whose cohomology objects are the Hodge structures $\H^i(X^\an, \Q(m))$. 
\lemm \mylabel{lemm_Deligneweak}
For $X / \Q$ smooth and projective and any $i, m$, we have 
\eqn \mylabel{eqn_Deligneweak}
\Hw^i (\underline \RG(X, m)) = \HD^i (X, m).
\xeqn
The induced isomorphism $\det \Hw^*(\underline \RG(X, m)) = \det \HD^* (X, m)$ respects the $\Q$-structure.  
\xlemm
\pf
The Hodge structures $L_i := \H^i (\underline{\RG}(X, m)) = \H^i(X^\an, \Q(m))$ are pure of weight $i-2 m$. For $i - 2m < 0$, $\Hw^0(L_i) = \Gamma_\MHS(L_i) = 0$. By duality, $\Hw^1 (L_i) = \Hw^0(L_i\dual(1))\dual = 0$ for $i - 2m >-2$. Hence, by \refeq{specweak}, 
$$\Hw^i \underline{\RG}(X, m) = \left \{ \begin{array}{ll} \Hw^1 (L_{i-1}) & i - 2m < 0 \\
\Hw^0 (L_i) & i-2m \geq 0 \end{array} \right.$$
The map in the exact sequences \refeq{longDeligne} between Betti and de Rham cohomology is the one from the definition of $\RGw(L_*)$. This shows \refeq{Deligneweak}. The identification of the $\Q$-structures follows similarly.
\xpf

For archimedean factors of $L$-functions of arbitrary motives, we use the Hodge realization functor (see \cite[Section 3]{Beilinson:Notes} for an early avatar):
\eqn \mylabel{eqn_Hodgerealization}
\RG_\mathrm{H}: \DM_{\Beilinson, c}(\Q)^\op \stackrel {\cong, \text{\refeq{DMgm}}} \lr \DMgm(\Q)^\op \r \DbH.
\xeqn
The right hand functor is Huber's Hodge realization functor \cite[2.3.5]{Huber:Realization}. It maps $\Mgm(X)(-m)$ to $\underline{\RG} (X, m)$. For any $M \in \DM_{\Beilinson, c}(\Q)$, the natural map $\RG_\mathrm H (\IHom(M, \one)) \r \IHom(\RG_\mathrm H (M), \RG_\mathrm H (\one))$ is an isomorphism. It is enough to check this on generators $M=\M(X)$ with $X/\Q$ smooth and projective, where it follows from $(\M(X))\dual = \M(X)\twi{\dim X}$. We obtain $\RG_\mathrm{H}(M\dual(1)) = (\RG_\mathrm{H} (M))\dual(-1)$. We put
\eqn \mylabel{eqn_RGw}
\RG_{\mathrm{wH}} := \RGw \circ \RG_\mathrm{H} : \DM_{\Beilinson, c}(\Q) \r \DbRQdet.
\xeqn
The composition of these functors with $\eta^*: \DM_{\Beilinson, c}(\Z) \r \DM_{\Beilinson, c}(\Q)$ will be denoted the same.

Finally, we recall the construction of the Deligne cohomology spectrum $\HD$ \cite{HolmstromScholbach:Arakelov}. We also sketch how to obtain similar spectra for Betti and de Rham cohomology. The aim is \refeq{DeligneQ}, the $\Q$-structure on Deligne cohomology groups of general motives.   

Let $\cC$ be either the category $\Sm^{G, \an}$ or $\Sm / \Q$. Consider simplicial presheaves $C(p)$ of pointed sets on $\cC$, for each $p \geq 0$, together with a ``product'' map $\cdot_C : C(p) \wedge C(p') \r C(p+p')$. Moreover, we assume there is an element $c_1 \in C(1)(\Gm)$, that restricts to zero at the point $1 \in \Gm$ (equivalently, a pointed map $c_1: (\Gm, 1) \r C(1)$) such that for any two maps $f_i: U \r \Gm$, $U \in \cC$, $i = 1, 2$,
\eqn \mylabel{eqn_c1}
f_1^* (c_1) \cdot_C (f_2^*(c_1) \cdot_C c') = f_2^* (c_1) \cdot_C (f_1^*(c_1) \cdot_C c').
\xeqn
The element $c_1$ is referred to as a \Def{bonding element}. Under these assumptions, the presheaves $C(p)$ with the bonding maps $\Gm \wedge C(p) \stackrel{c_1 \wedge \id} \lr C(1) \wedge C(p) \stackrel{\cdot_C}\lr C(p+1)$ form a symmetric $\Gm$-spectrum $C$ (where the $\Sigma_p$-action on $C(p)$ is trivial). The category of such spectra is denoted $\Spt(\cC)$. It is endowed with a model structure whose homotopy category $\SH(\Q)$ (or $\SH(\RR^\an)$) satisfies (cf.\ e.g.\ \cite[Section 1]{Ayoub:Note} for the analytic version):  
$$\Hom_{\SH}(\Sigma^\infty (X \sqcup \{*\}) \wedge S^n \wedge \Gm^{\wedge m}, C) = \pi_{n+m+N}(C(m+N)(X))$$
for any $X \in \cC$, and $n, m \in \Z$ and $N \gg 0$, provided that  
\begin{enumerate}
\item all levels $C(p)$ are homotopy invariant: $C(p)(-) \r C(p)(- \x \A)$ (respectively, $- \x (\A)^\an$) is a weak equivalence,
\item all levels $C(p)$ satisfy descent (with respect to the Nisnevich and the analytic topology, respectively), and
\item $C$ is an $\Omega$-spectrum. In the presence of the first two conditions, this is implied by the bundle formula, which says that
$$\oplus_{i=0}^1 p_X^*(-) \cdot_C p_{\Gm}^* (c_1)^i : \oplus \pi_{i+*}(C(p-i)(X)) \r \pi_{*} (C(p)(X \x \Gm))$$
is an isomorphism, where $p_X$, $p_{\Gm} : X \x \Gm \r X$, $\Gm$ are the projections. 
\end{enumerate}
The spectra below are all obtained by putting $C(p) := DK (\tau_{\geq 0} A(p))$ for appropriate complexes of abelian groups $A(p)$. Here $\tau$ is the good truncation and $DK$ the Dold-Kan equivalence.

We now define four different (but isomorphic) spectra representing Betti cohomology with real coefficients by specifying the levels $C(p)$ and the bonding elements in $C(1)(\Gm)$. The product structure map on the level complexes is obvious for these Betti cohomology spectra, and is strictly commutative and associative. For any presheaf of abelian groups $F$ on $\Sm^{G, \an}$, we define the \v Cech-complex in degrees $n \geq 0$  
$$\mathfrak C^n F (X) := \varprojlim F (U^{n+1}).$$
The limit runs over the directed category of all open covers $\{U_i \}$ of $X \in \Sm^{G, \an}$ and $U := \sqcup U_i$. 
Given some involution $\ol{?} : F \r F$, we write $\mathfrak C^G F$ for the subcomplex consisting of elements that are fixed by $\ol{\Fr_\infty^*}$. 

Let $\HBR^{(1)}$ be the spectrum whose levels are $\mathfrak C^G (\RR(p)[p])$. To describe the bonding element, we replace $\Gm^\an$ by $S^1$ (equipped with its usual topology). The inclusion $S^1 \subset \Gm^\an$ is a homotopy equivalence, and an explicit description of a \v Cech cocycle generating $\H^1(\Gm^\an, \CC)$ is left to the reader. As for $S^1$, consider the standard covering by $U_{\pm} = \{z \in S^1, \pm \Re(z) > -0.5 \}$. This covering is equivariant with respect to $z \mapsto \ol z$. Frobenius $\Fr_\infty$ acts on the \v Cech complex
$$\RR(1)(U_+) \oplus \RR(1)(U_-) \r \RR(1)(U_+ \cap U_-) = \RR(1)^2, (a, b) \mapsto (v,w) := (b-a, b-a)$$
as $(a, b) \mapsto (a,b)$ and $(v,w) \mapsto (w,v)$. Hence $(\pi i, -\pi i) \in \RR(1)(U_+ \cap U_-)$ is a $\ol{\Fr_\infty}$-invariant element which generates $\H^1(\Gm^\an, \RR(1))^G$. This determines the spectrum $\HBR^{(1)}$. It is well-known that $\H^* (\mathfrak C^*(\RR)(X))=\H^*(X, \RR)$. Thus
\eqn \mylabel{eqn_HBR1}
\Hom_{\SH(\RR^\an)} (\Sigma^\infty X, \HBR^{(1)}(p)[n]) = \HBetti^n(X, \RR(p))^{(-1)^p},
\xeqn
where the superscript at the right denotes the subgroup of elements $a$ satisfying $\Fr_\infty^*(a)=(-1)^p a$. The complexes $\Tot (\mathfrak C^G (E_\RR^* (p)[p]))$ and the bonding element induced by the previous one via the inclusion $\RR(1)[1] \subset E_\RR^0(1)[1]$ yield a spectrum $\HBR^{(2)}$ that is naturally isomorphic to $\HBR^{(1)}$, since $\RR \r E_\RR^*$ is a quasi-isomorphism of sheaf complexes. Consider the spectrum $\HBR^{(3)}$ whose levels are the one of $\HBR^{(2)}$, but the bonding element is the $1$-form 
$$dz/z \in E^1_{\RR, \P}(\log \{0, \infty\}) \r \mathfrak C^0 E_\RR^1(1)(\Gm^\an) \subset \Tot(\mathfrak C^* E_\RR^*(1))^1(\Gm^\an).$$ 
Both $\HBR^{(2)}$ and $\HBR^{(3)}$ are $\Omega$-spectra (the above bonding element and $dz/z$ give the same element in $\H^1(\Gm^\an, \RR(1))$ by Cauchy's residue formula). The identity map between their level-$0$-complexes thus yields a canonical isomorphism of spectra (in $\SH(\RR^\an)$). 
The complexes $E_\RR^{*,G}(p)[p]$ (again $?^G$ denotes invariants under $\ol{\Fr_\infty^*}$) together with the bonding element $dz/z$ form a spectrum denoted $\HBR^{(4)}$. The obvious quasi-isomorphism $E_\RR^{*, G} = \mathfrak C^{0, G} E_\RR^* \r \Tot (\mathfrak C^G E^*_\RR)$ induces an isomorphism $\HBR^{(4)} \r \HBR^{(3)}$ in $\SH(\RR^\an)$. The purpose of the chain of isomorphisms $\HBR^{(4)} \cong \HBR^{(1)}$ is the existence of $\HBQ^{(1)}$, the obvious $\Q$-linear variant of $\HBR^{(1)}$. It induces a $\Q$-structure on the groups represented by $\HBR^{(4)}$. 

As for de Rham cohomology, consider the complexes $E^F(p)^G := \cone (F^p E^* \r E^*)^G[p-1]$. The product
\eqnarr
(F^p E^n \oplus E^{n-1}) \t (F^{p'} E^{n'} \oplus E^{n'-1}) & \r  & (F^{p+p'} E^{n+n'} \oplus E^{n+n'-1}), \\
(f_1, e_1) \t (f_2, e_2) & \mapsto & (f_1 \wedge f_2, f_1 \wedge e_2)
\xeqnarr 
is strictly associative, but in general commutative only up to homotopy \cite[Section 3]{EsnaultViehweg}. However, putting $c_1 = (dz/z, 0) \in E^F(1)(\Gm) = (F^1 E^1 \oplus E^0)(\Gm)$, \refeq{c1} clearly holds. We obtain a spectrum $\HdR^{/F, \an} \in \SH(\RR^\an)$. Using $\Fr_\infty^*$-invariant algebraic differential forms, i.e., $\Omega^{*, \alg, G}$ instead of $E^{*, G}$, we get a similar spectrum denoted $\HdR^{/F, \alg} \in \SH(\Q)$. For smooth $X/\Q$, the obvious maps 
$$\Omega^{*, \alg}(X) \t_\Q \RR \r \Omega^{*, \alg}(X_\RR) = \Omega^{*, \alg, G}(X_\CC) \leftarrow \Omega^{*, G}(X) \r E^{*, G}(X)$$ 
are filtered (with respect to the Hodge filtration) quasi-isomorphisms by flat base change for $\Omega^{*, \alg}$, GAGA and \cite[Thm. 2.1]{Burgos:Dolbeault}. We thus get an isomorphism $c_* \an_* \HdR^{/F, \an} = \HdR^{/F, \alg} \t \RR$ in $\SH(\Q)$. Here $c: \Spec \RR \r \Spec \Q$.  
      
Finally, the complex
\eqn \mylabel{eqn_Delignep}
\Deligne(p)^G := \cone (E_\RR^{*,G} (p)[p] \r E^F(p)^G)[-1]
\xeqn
carries a product map $\cdot_{\Deligne, \alpha}$ depending on some auxiliary parameter $\alpha \in \RR$. It is only commutative and associative up to homotopy (for each $\alpha$). Again, $c_1 = (dz/z, dz/z, 0) \in \Deligne(1)^0(\Gm) = (E_\RR^{n+p}(p) \oplus F^p E^{n+p} \oplus E^{n+p-1})(\Gm)$ satisfies \refeq{c1} (independently of $\alpha$, see the multiplication table in \lc). The resulting spectrum $\HD$ sits in a distinguished triangle in $\SH(\RR^\an)$, $\HD \r \HBR^{(4)} \r \HdR^{/F, \an}$ and thus, in $\SH(\Q)$,
\eqn \mylabel{eqn_distinguishedHD}
c_* \an_* \HD \r c_* \an_* (\HBQ^{(1)} \t_\Q \RR) \r \HdR^{/F, \alg} \t_\Q \RR \r c_* \an_* \HD[1].
\xeqn  
From now on, we write $\HD$ for $\an_* c_* \HD \in \SH(\Q)$. This is the spectrum established in \cite[Section 3]{HolmstromScholbach:Arakelov}, except for two inessential differences: instead of $\Deligne(p)$, \lcs used other complexes that are homotopic (including the product structure, regardless of $\alpha$) to $\Deligne(*)$. Secondly, the construction of \lcs builds a symmetric $\P$-spectrum, but again this is inessential at the level of the homotopy categories, since $- \wedge \P = - \wedge \Gm \wedge S^1_s$, where $S^1_s$ is the simplicial sphere. By \cite[Thm. 3.6]{HolmstromScholbach:Arakelov}, 
$$\Hom_{\DM_{\Beilinson}(\Q)} (\M(X), \HD(p)[n]) = \Hom_{\SH(\Q)}(\M(X), \HD(p)[n]) = \HD^n(X, p)$$
for any $X \in \Sm/\Q$. For any $M \in \DM_{\Beilinson, c}(\Q)$, \refeq{distinguishedHD} induces an isomorphism 
\eqnarra \mylabel{eqn_DeligneQ}
\det \HD^* (M) & = & \det c_* \an_* \HBR^{(4), *} (M) \t \det^{-1} c_* \an_* \HdR^{/F,\an} (M) \nonumber \\
& = & \left (\det c_* \an_* \HBQ^{(1), *} (M) \t \det^{-1} \HdR^{/F, \alg, *}(M) \right) \t_\Q \RR. 
\xeqnarra
Here $\det \HD^*(M) := \t_{n \in \Z} \det^{(-1)^n} \Hom(M, \HD[n])$ etc. is well-defined since $M$ is compact. This is the promised extension of \refeq{Qstructure} to Deligne cohomology groups of general geometric motives.  

\mylabel{exactnessBetti}Applied to the Betti realization, the exactness axiom (see \refax{mixed}\refit{exactness}) means 
\eqn \mylabel{eqn_exactnessBetti}
\Hom(M, \HBR)=\Hom(\pH^0(M), \HBR), \ \ \text{ for all }M \in \DMBeic(\Q)
\xeqn 
and likewise for de Rham cohomology. This implies that for any smooth projective $X_\eta / \Q$,
\eqn \mylabel{eqn_Delignemixed}
\HD^i(\h^{-b-1}(X_\eta, -m)) = \begin{cases}
\HD^{b+1}(X_\eta, m) & i = 0 \text{ and } b+1-2m \geq 0 \\
\HD^{b+2}(X_\eta, m) & i = 1 \text{ and } b+1-2m \leq -2\\
0 & \text{else.}
\end{cases}
\xeqn

\subsection{Arakelov motivic cohomology}
\mylabel{sect_Arakelov}

In order to formulate \refco{L} below, we need to recall some facts about Arakelov motivic cohomology. 
\theo \cite{HolmstromScholbach:Arakelov, Scholbach:ArakelovII} \mylabel{theo_Arakelovmot}
In $\DM_{\Beilinson}(\Z)$, there is a unique map $\ch: \one \r \eta_* \HD$ representing the Chern class map from motivic cohomology to Deligne cohomology, i.e. 
$$\Hom_{\DM_{\Beilinson}(\Z)}(\M(X), \one(p)[n]) \stackrel {\ch(p)[n]}\lr \Hom_{\DM_{\Beilinson}(\Z)}(\M(X), \eta_* \HD(p)[n])$$ 
agrees with the Chern class $K_{2p-n}(X)^{(p)}_\Q \r \HD^n(X_\Q, p)$ (also known as Beilinson regulator) for all regular projective schemes $X/\Z$. There is a certain, canonically defined object $\onehat \in \DM_{\Beilinson}(\Z)$ called \emph{Arakelov motivic cohomology spectrum} such that there is a distinguished triangle   
\eqn \mylabel{eqn_onehat}
\onehat \stackrel{f} \r \one \stackrel \ch \lr \eta_* \HD \stackrel \delta \lr \onehat[1].
\xeqn
Moreover, given another object $\onehat'$ in a similar triangle, there is a \emph{unique} isomorphism $\onehat \r \onehat'$ in $\DM_{\Beilinson}(\Z)$ fitting in the obvious commutative diagram of distinguished triangles.  
\xtheo

\defi
Given a motive $M \in \DMBeic(\Z)$, its \Def{Arakelov motivic cohomology} is defined as
$$\Hhat^i(M, m) := \Hom_{\DM_{\Beilinson}(\Z)} (M, \onehat(m)[i]).$$
We write $\Hhat^i(X, m) := \Hhat^i(\M(X), m)$. We also consider the $\RR$-linear variant of these groups, denoted $\Hhat_\RR^i(X, m)$, obtained by replacing $\one$ by $\one_\RR$ in \refeq{onehat}. This amounts to tensoring the motivic cohomology groups with $\RR$.
\xdefi 

The triangle \refeq{onehat} induces long exact sequences 
\eqn \mylabel{eqn_Hhat}
\Hhat^i_\RR(M, m) \r \H^i(M, m)_\RR \r \HD^i(M, m) \r \Hhat_\RR^{i+1}(M, m).
\xeqn

On the other hand, we have the notion of arithmetic $K$-theory. For a regular and projective scheme $X$ over $\Z$, such groups $\Khat^T_n(X)$ have been defined by Gillet and Soul\'e for $n=0$ and for higher $n$ by Takeda \cite[Section 6]{GilletSoule:CharacteristicII}, \cite[p. 621]{Takeda:Higher}. These groups sit in an exact sequence
$$K_{n+1}(X) \r \oplus_{p \in \Z} \Deligne(p)^{2p-n-1, G}(X) / \im d_\Deligne \r \Khat^T_n(X) \r K_n(X) \r 0$$
where $\Deligne(p)^G$ is the complex defined in \refeq{Delignep}. Moreover, they come with a Chern class map $\ch: \Khat^T_n(X) \r \oplus_{p \in \Z} \Deligne(p)^{2p-n, G}(X)$. The group $\Khat_n(X) := \ker \ch$ fits in a long exact sequence
\eqn
\mylabel{eqn_Khat}
\dots \r \oplus_{p \in \Z} \HD^{2p-n-1}(X, p) \r \Khat_n(X) \r K_n(X) \r \oplus \HD^{2p-n}(X, p) \r \dots
\xeqn
The group $\Khat_0^T(X)_\Q$ is also isomorphic, via the arithmetic Chern class to $\oplus_p \CHhatGS^p(X)_\Q$, where $\CHhatGS$ denotes the arithmetic Chow group of Gillet and Soul\'e \cite[3.3.4]{GilletSoule:Arithmetic}. It is generated by arithmetic cycles $(Z, g_Z)$, i.e., cycles $Z \subset X$ and Green currents, i.e., such that $\omega_Z := \delta_Z - 2 \partial \ol \partial g_Z$ is a differential form. Here $\delta_Z$ is the Dirac current.   
Under the arithmetic Chern class, the subgroup $\Khat_0(X)_\Q \subset \Khat_0^T(X)_\Q$ corresponds to the subgroup $\CHhat^*(X) \subset \CHhatGS^*(X)$ generated by arithmetic cycles $(Z, g_Z)$ such that $\omega_Z =0$ \cite[Thm. 7.3.4]{GilletSoule:CharacteristicII}.    

For a smooth scheme $X$ over $S \subset \Spec \Z$, the resulting decomposition of $\Khat_0(X)_\Q$ in Adams eigenspaces is extended to higher $\Khat$-theory \cite[Cor. 6.2]{Scholbach:ArakelovII}: $\Khat_n(X)_\Q$ decomposes as a direct sum of Adams eigenspaces $\oplus \Khat_n(X)^{(p)}_\Q$, compatibly with \refeq{Khat}. In fact, this statement is derived from a canonical isomorphism 
\eqn \mylabel{eqn_KhatCHhat}
\Hhat^i(X, m) = \Khat_{2m-i}(X)^{(m)}_\Q \ \ (= \CHhat^m(X)_\Q \text{ for }i=2m).
\xeqn
 
\defi \mylabel{defi_A.i.p.}
Let $S \subset \SpecZ$ be an open subscheme and let $M \in \DM_{\Beilinson}(S)$ be any motive. The natural pairing of motivic homology (see \refeq{motivichomology}) and Arakelov motivic cohomology, 
$$\pi_M: \H_{-2}(M, -1)_\RR \x \Hhat^{0}_\RR(M) \r \Hhat^2_\RR(\one_S, 1)$$
given by the composition of morphisms in $\DM_{\Beilinson}(S)$ is called \Def{Arakelov intersection pairing}.
\xdefi

\rema 
\mylabel{rema_Hhat.basic}
\begin{enumerate}[(i)]
\item
For $M \in \DMBeic(S)$, we often tacitly identify $\H_{-2}(M, -1) \cong \H^2(M\dual, 1)$, cf.\ \refeq{reflexive}.
\item 
The Arakelov intersection pairing is functorial in $M$ in an obvious sense. 
\item \mylabel{item_duality}
Let $M \in \DMBeic(S)$. Consider
\eqn \mylabel{eqn_dualitycompatible}
\begin{array}{cccccl}
\Hhat_\RR^{0}(M) & \x & \H_\RR^{2}(M\dual, 1) & \lr & \Hhat^2(\one, 1)  \\
\downarrow & & \uparrow & & \downarrow = \\
\H^0_\RR(M) & \x & \Hhat_\RR^{2}(M\dual, 1) & \lr & \Hhat^2(\one, 1)  \\
\downarrow & & \uparrow & & \uparrow \cong \\
\HD^0(M) & \x & \HD^{1}(M\dual, 1) & \lr & \HD^1(\one, 1), \end{array}
\xeqn
where in the first row $(a: M \r \onehat, b: M\dual \r \one\twi{-1})$ is mapped to $\mu \circ (a \t b) \circ \text{coev}$, where the coevaluation $\one \r M \t M\dual$ is obtained from \refeq{reflexive}, $\mu: \one \t \onehat \r \onehat$ is the $\one$-module structure map for $\onehat$. This is just another way to write $\pi_M$. Likewise, the second row pairing is $\pi_{M\dual\twi{-1}}$. The pairing in the third row is defined similarly using the product of the ring spectrum $\mu_\Deligne: \HD \t \HD \r \HD$ instead. This diagram is commutative. This follows from the commutativity of the following diagram, which in turn is a rephrasing of the fact that \refeq{onehat} is a distinguished triangle of $\one$-modules. 
$$
\xymatrix{
\onehat \t \onehat \ar[d]^{\id \t f} \ar[r]^{f \t \id} &
\one \t \onehat \ar[d]^\mu &
 &
\one \t \HD[-1] \ar[ll]^\delta \ar[d]^{\ch \t \id} \\
\onehat \t \one \ar[r]^\mu &
\onehat &
\HD[-1] \ar[l]^\delta & 
\HD \t \HD[-1] \ar[l]^{\mu_\Deligne}. 
} 
$$
\item \mylabel{item_duality2}
The pairing $\HD^0(M) \x \HD^1(M\dual, 1) \r \RR$ is a perfect pairing for any $M$. It suffices to see this for $M = \M(X)(p)[n]$ for $X / \Z$ regular and projective, in which case it follows from the identification of Deligne cohomology with weak Hodge cohomology (\refle{Deligneweak}) and the duality of weak Hodge cohomology, \refeq{dualityweak}.
This plays an important role in the compatibility of our $L$-values conjecture with respect to the functional equation, see \refth{structure}\refit{2}. 
\end{enumerate}
\xrema

\exam \mylabel{exam_Fp}
Consider a motive $M = i_* N$, where $i: \SpecFp \r \SpecZ$ and $N \in \DMBeic(\Fp)$ (for example $M = \M(\Fp) = i_* i^* \one \twi{-1}$). The forgetful map $f: \Hhat^0_\RR(M) \r \H^0(M)_\RR = \H^{-2}(N, -1)_\RR$ is an isomorphism and the pairing $\pi_M$ coincides with the natural pairing $\H_{-2}(N, -1)_\RR \x \H^{-2}(N, -1)_\RR \r \H^0(\one_{\Fp}, 0)_\RR = \RR$ followed by the pushforward $i_* : \Hhat^0(\one_{\Fp}, 0)_\RR \r \Hhat^2_\RR(\one_\Z, 1)$, which is $\log p: \RR \r \RR$ \cite[Theorem 6.4.(i)]{Scholbach:ArakelovII}.
\xexam 

\exam \mylabel{exam_AipX}
Let $X$ be a regular projective scheme over $S \subset \SpecZ$ of constant dimension $d$. We pick some open $j: U \subset S$ such that $X_U$ is smooth over $U$. Let $M := \M(X)\twi{m-d}[i] \in \DMBeic(S)$. Then 
$\H_{-2}(M, -1) = 
K_i(X)^{(m)}_\Q$ by absolute purity. Let $M_U := j^* M \in \DMBeic(U)$. Consider 
$$\xymatrix{
K_{i}(X_U)^{(m)}_\Q \x \Khat_{-i}(X_U)^{(d-m)}_\Q \ar[r]^(.7)\cup & \Khat_0(X_U)^{(d)}_\Q \ar[r]^{{f_U}_*} & \Khat_0(U)^{(1)}_\Q \\
\H_{-2}(M_U, -1) \x \Hhat^{0}(M_U) \ar[u]^\cong \ar[rr]^{\pi_{M_U[-i]}} && \Hhat^2(\one_U, 1) = \RR / \sum_{p \notin U} \log p \Q \ar[u]^\cong \\
\H_{-2}(M, -1) \ar@{->>}[u]^{j^*} \x \Hhat^{0}(M) \ar@{->>}[u] \ar[rr]^{\pi_{M[-i]}} & & \Hhat^2(\one_S, 1) = \RR / \sum_{p \notin S} \log p \Q \ar@{->>}[u]^{j^*} 
}$$
In the first row, the pushforward ${f_U}_*$ is not the pushforward on arithmetic $K$-theory, but the one on arithmetic Chow groups using the arithmetic Chern class isomorphism \refeq{KhatCHhat}. The top square is commutative by \cite[Thm. 7.4.]{Scholbach:ArakelovII}. The bottom square is commutative by definition. See also \refre{height}.

\xexam

\section{$L$-functions of motives over number rings} 
\mylabel{sect_LOF}
Let $F$ be a number field and $\OF$ its ring of integers. For every finite prime $\pp$ of $\OF$ we fix a rational prime $\ell$ that does not lie under $\pp$. Moreover, fix for every $\ell$ an embedding $\sigma_\ell : \Ql \r \CC$. All subsequent definitions of $L$-functions are taken with respect to these choices.

\defi \mylabel{defi_LfunctionoverOF}
The $L$-series of a mixed motive $M_\eta$ over $F$ is defined by
$$ L_F(M_\eta, s) := \prod_{\pp < \infty} \det \left( \Id - \Fr^{-1} \cdot N(\pp)^{- s} | ({M_\eta}_\ell \t_{\Ql, \sigma_\ell} \CC)^{I_\pp} \right)^{-1}.$$
The $L$-series of a geometric motive $M$ over $\OF$ is given by
$$L_{\SpecOF}(M, s) := L(M, s) := \prod_{\pp < \infty} \det \left(\Id - \Fr^{-1} \cdot N(\pp)^{-s}|(i_\pp^! M)_\ell \t_{\Ql, \sigma_\ell} \CC \right)^{-1}.$$
The first definition is classical, the second is a natural adaptation to motives over $\OF$. The products run over all finite primes of $\OF$, $\Fr$ is the arithmetic Frobenius map (given on residue fields by $a \mapsto a^{N(\pp)}$), $N(\pp)$ is the norm of $\pp$, $i_\pp$ denotes the immersion of the corresponding closed point and $-_\ell$ denotes the $\ell$-adic realization functor, see \refeq{ell-adic}. The determinants are understood in the sense of \refsect{LA}. 
The superscript $I_\pp$ denotes the invariants under the action of the inertia group.
\xdefi

\bem \mylabel{bem_propertiesL}
By \refaxtstructureRealizations, the $\ell$-adic realization ${M_\eta}_\ell$ is in fact an $\ell$-adic sheaf. For example, $(\h^{-b-1}(X_\eta, -m))_\ell = \H^{b+1}(X_\eta, \Ql(m))$ for some scheme $X_\eta$ over $F$. 

The independence of the choices of $\ell$ and the embeddings $\sigma_\ell$ is discussed around \refle{independencel}. See also \refth{reductionperfFp}.

The $L$-series for motives over $\OF$ is multiplicative, i.e., given a triangle $M \r M' \r M''$ in $\DMBeic(\OF)$, one gets 
$$L(M', s) = L(M, s) \cdot L(M'', s).$$ 
A similar property does \emph{not} hold for $L$-functions of motives over $F$ \cite{Scholl:Remarks}.
\iffinal
See also \cite[1.3.3]{FPR}.
\fi   
\ifpreprint 
Scholl's notion of mixed motives over $\Z$ (which are not the same thing as mixed motives over $\Z$ in our sense, but mixed motives over $\Q$ with certain additional non-ramification properties), as well as Fontaine's and Perrin-Riou's notion of $f$-exact sequences \cite[1.3.3]{FPR} are designed to grapple with this phenomenon. 
\fi

By definition and the calculation of $\ell$-adic cohomology of $\P_\Fpp$
\ifpreprint
\cite[Example VI.5.6 and p. 163]{Milne:EC}, 
\fi
\iffinal
, 
\fi
one has 
\eqn \mylabel{eqn_Lfunctionshift}
L(M(-m), s) = L(M, m+s), \ \ m \in \Z .
\xeqn

For an open subscheme $j : \SpecOF \backslash Z \r \SpecOF$ with complement $i: Z \r \SpecOF$, the $L$-function of $j_* j^* M$ is the one of $M$, but the Euler factors for the points in $Z$ are omitted. This follows from $i^! j_* = 0$.
\xbem

The following lemma is well-known, see \cite[Prop. 3.8.(ii)]{Deligne:Constantes} or \cite[VII.10.4.(iv)]{Neukirch:AZT} for similar statements. It permits to replace any number ring $\OF$ by $\Z$ and to study $L$-values of motives over $\Z$, only. 

\lemm \mylabel{lemm_pushdown}
The $L$-series is an absolute invariant of a motive, i.e., for any geometric motive $M$ over $\SpecOF$ we have $L_{\SpecOF}(M, s) = L_{\SpecZ}(f_* M, s)$, where $f : \SpecOF \r \SpecZ$ denotes the structural map.
\xlemm
\ifpreprint
\pf
By definition, one can reduce to the case where $M$ is supported on a single prime $\pp$ in $\OF$, that is, $M = {i_\pp}_* M'$. Using $(f_* M)_\ell = f_! (M_\ell) = f_* (M_\ell)$, one reduces the claim to a question about representations of the Galois groups of the involved finite fields. Then the claim follows from a linear-algebraic calculation done e.g.\ in the proof of \cite[VII.10.4.(iv)]{Neukirch:AZT}. The details are omitted.
\xpf
\fi

We now relate $L$-series of motives over $\Q$ to ones over $\Z$. Recall the notion of smooth motives from \refde{generically.smooth}. The following lemma is proven in \cite[\fcohoSectionIntermediateExtensionEll]{Scholbach:fcoho} as a corollary of the exactness axiom for $\ell$-realization functors (see around \refeq{exactnessell}). 

\lemm  \mylabel{lemm_excreal}
Let $M$ be a mixed smooth motive over $U$, where $j : U \r \SpecZ[1/\ell]$ is an open subscheme. Let $i$ be the complementary closed immersion to $j$ and let $\eta'$ and $\eta$ be the generic point of $U$ and $\SpecZ[1/\ell]$, respectively. Then $(i^! j_{!*} M)_\ell = i^* (\R^0 \eta_* \eta'^* M_\ell[1])[-1]$.
\xlemm

The following proposition relates $L$-series of motives over $\Q$ and $\Z$. 
Our main example is $M_\eta = \h^{-b-1}(X_\eta, -m)$ and $M=\h^{-b}(X, -m)$ where $X / \Z$ is some projective scheme whose generic fiber $X_\eta / \Q$ is smooth.

\satz \mylabel{satz_corrj}
Let $M_\eta \in \MM(\Q)$. Pick some $M \in \MM(\Z)$ with $M_\eta = \eta^*[-1] M$. Then
$$L_\Q(M_\eta, s)^{-1} = L_\Z(\eta_{!*} \eta^* M, s).$$
\xsatz

\pf
For sufficiently small $j : U \r \SpecZ$, the right hand side is equal to
\eqnarr
L_\Z(j_{!*} j^* M, s) 
& \stackrel{\text{\ref{lemm_excreal}}}= &  \left (\prod_p \det \left (\Id - \Fr^{-1} p^{-s}|i_p^* (\R^0 \eta_* \eta^* M_\ell [1])[-1] \right) \right )^{-1} \\
& = &  \prod_p \det \left (\Id - \Fr^{-1} p^{-s}|i_p^* \R^0 \eta_* {M_\eta}_\ell \right) \\
& = &  \prod_p \det \left (\Id - \Fr^{-1} p^{-s}|({M_\eta}_\ell)^{I_p} \right) 
 = 
L_\Q(M_\eta, s)^{-1}.
\xeqnarr
\xpf

\subsection{Hasse-Weil $\zeta$-functions -- Motives with compact support}

\mylabel{sect_zetafunctions}

\defi (see e.g.\ \cite{Serre:Zeta})
The \Def{Hasse-Weil zeta function} of a scheme $X / \Z$ (always separated and of finite type) is defined as $\zeta(X, s) := \prod_x (1 - N(x)^{-s})^{-1}$. The product is over all closed points $x$ of $X$, and $N(x)$ denotes the cardinality of the (finite) residue field of $x$.
\xdefi

Recall from \refeq{motiveX} the motive with compact support $\Mc(X)$ of some scheme $X$.

\satz \mylabel{satz_HasseWeil}
For any scheme $X / \Z$, we have
$$\zeta(X, s) = L(\Mc(X), s).$$
\xsatz
\pf
Writing $i_p: \SpecFp \r \SpecZ$ and $X_p := X \x \Fp$, base-change implies $i_p{}_* i_p^! \Mc(X) = \Mc(X_p)$. (At the right hand side, $X_p$ is seen as a $\Z$-scheme.) Therefore, $L(\Mc(X), s) = \prod_p L(\Mc(X_p), s)$. A similar decomposition for the $\zeta$-function allows us to assume that $X$ is an $\Fp$-scheme. The $\ell$-adic realization functor satisfies $(f_* f^! \one)_\ell = f_! f^* \Ql$. Grothendieck's trace formula (see e.g. \cite[Sections VI.12, 13]{Milne:EC}) says 
\eqnarra
\zeta(X, s) & = & \prod_{i=0}^{2 \dim X} \left(\det \left(\Id - \Fr^{-1} \cdot p^{-s}|\H^i_{c}(X \x_{\Fp} \Fpq, \Ql) \right)\right)^{(-1)^{i+1}} \nonumber \\
& = & \det (\Id - \Fr^{-1} \cdot p^{-s} | f_! f^* \Ql)^{-1}, \nonumber 
\xeqnarra
where $\H^i_c(X \x \Fpq, \Ql) = \H^i (f_! f^* \Ql)$ denotes $\ell$-adic cohomology with compact support. 
\xpf



The $L$-series of a motive over $\Q$ is conjectured to be independent of the choice of $\ell$ and $\sigma_\ell$ in every factor (assuming $p \neq \ell$). This is known for the individual Euler factors at $p$ if the motive is $\h^i(X_\eta, n)$, where $X_\eta$ is a variety with good reduction at $p$, by Deligne's work on the Weil conjectures \cite[Th. 1.6]{Deligne:Weil1}. From \refsa{HasseWeil} we now immediately obtain another statement concerning independence of $\ell$.

\defi \mylabel{defi_accessible}
The smallest triangulated subcategory of $\DMBeic(\Z)$ containing the motives $\M(X)(n)$ ($n \in \Z$) of all regular schemes $X$ which are projective and flat over $\Z$, and the image of $i_* : \DMBeic(\Fp) \r \DMBeic(\Z)$ for all primes $p$, is called
$\DMtrgm(\Z)$\index{DMgm2@$\DMtrgm$} and called category of \emph{accessible motives}.
Its triangulated subcategory generated by $\M(X)(n)$ where $X$ is regular and projective, but not necessarily flat over $\Z$ (such as a smooth projective $X / \Fp$) is called the category of \emph{easily accessible motives}. 
\xdefi 

\rema \mylabel{rema_truly.geometric}
\begin{enumerate}[(i)]
\item
By de Jong's resolution of singularities using alterations, the 
thick closure (i.e., closure under direct summands and triangles) of the category of easily accessible motives contains the motives $\M(X)(n)$ of all $X$ schemes (of finite type) over $\Z$. Therefore, this thick closure is the entire category $\DMBeic(\Z)$ of geometric motives.
\item \mylabel{item_truly.geometric.generators}
By the proof of \cite[\fcohoSatzgeneratorsoverOF]{Scholbach:fcoho}, $\DMtrgm(\Z)$ is contained in the triangulated category generated by 
$i_* \DMBeic(\Fp)$ and motives of the form $E := \eta_{!*} \eta^* \h^{-b}(X, -m)$, where $X / \Z$ is regular, flat and projective. 
\end{enumerate}
\xrema

The following lemma shows that the question of independence of $L$-functions of $\ell$ is solely about the behavior of $L$-functions under direct summands. 

\lemm \mylabel{lemm_independencel}
For any easily accessible motive $M$ over $\Z$, the $L$-series $L(M, s)$ does not depend on the choices of $\ell$ (provided $p \nmid \ell$) and $\sigma_\ell$.
\xlemm

\pf
Using \refeq{Lfunctionshift}, we 
may assume $M = \M(X) = \Mc(X)$ for some $X$ which is projective over $\Z$ (and regular). Then the claim immediately follows from \refsa{HasseWeil}.
\xpf

\subsection{Meromorphic continuation and functional equation}\mylabel{sect_funceq}
Properties of $L$-series for motives over $\Q$ tend to generalize to ones over $\Z$, given that the property in question is known for motives over $\Fp$. We illustrate this by the absolute convergence, meromorphic continuation, and the functional equation. Recall from \cite[5.2.]{Deligne:Valeurs} or \cite[p.\ 4]{Schneider:Introduction} the definition of the archimedean Euler factor $L_\infty(V, s)$ for a mixed Hodge structure $V$. Essentially, $L_\infty(V, s)$ is a product of $\Gamma$-functions. The pole order at $s=0$ is given by (\cite[Lemma 7.1.]{Beilinson:Notes} or 
\cite[III.1.2.5 + III.1.2.3]{FPR}): 
\eqn
\mylabel{eqn_ordLinfty}
\ord_{s=0} L_\infty(V, s) = - \dim_\RR \Hw^1 (V\dual(1)).
\xeqn
For $V_* \in \D^\bound_\mathrm H$, we put $L_\infty(V_*, s ) := \prod_{i \in \Z} L_\infty(\H^i(V_*), s)^{(-1)^i}$. Here $\H^i(V_*)$ denotes the $i$-th cohomology Hodge structure of the complex $V_*$. 

\defi \mylabel{defi_completedL}
Let $M$ be a geometric motive over $\Z$ or a mixed motive over $\Q$. The function 
$$L_\infty(M, s) := L_\infty(\RG_\mathrm{H}(M), s)$$
is called the \Def{archimedean factor} of the $L$-function of $M$. Here $\RG_\mathrm{H}$ is the Hodge realization functor \refeq{Hodgerealization}. The \Def{completed $L$-function} of $M$ is defined as 
$$\Lambda(M, s) := L(M, s) L_\infty(M, s).$$
\xdefi

Much the same way as $L$-functions of motives over $\Q$, archimedean factors are not multiplicative with respect to short exact sequences of Hodge structures. (See \cite[1.1.9, 1.2.5]{FPR} for a necessary and sufficient criterion for multiplicativity.) 


The following is a long-standing conjecture concerning $L$-functions \cite{Deligne:Constantes}, \cite[5.2, 5.3]{Deligne:Valeurs} or \cite[p. 610, 699]{FPR}:
\conj \mylabel{conj_genpropLseries}
Let $M_\eta$ be a mixed motive over $\Q$. The $L$-series $L_\Q(M_\eta, s)$ converges absolutely for $\Re (s) \gg 0$ and has a meromorphic continuation to the complex plane. There is a functional equation relating the $\Lambda$-functions of $M_\eta$ and $M_\eta\dual(-1)$: 
$$\Lambda (M_\eta, s) = \epsilon (M, s) \Lambda (M_\eta\dual(-1), -s),$$
where $\epsilon(M, s)$ is of the form $a b^s$, with nonzero constants $a$ and $b$ depending on $M$.
\xconj

\lemm \mylabel{lemm_functionalequation}
\refco{genpropLseries} implies the following: for any accessible motive $M$ over $\Z$ (\refde{accessible}), the $L$-series $L(M, s)$ converges absolutely for $\Re(s) \gg 0$, has a meromorphic continuation to the complex plane, and there is a functional equation $\Lambda (M, s) = \epsilon (M, s) \Lambda (M\dual(-1), -s)$, where $\epsilon(M, s)$ is of the form $a b^s$, with nonzero constants $a$ and $b$ depending on $M$.
\xlemm
\pf
The  claim is triangulated, since the assignments $M \mapsto L(M, s)$, and $M \mapsto L_\infty(M, s) / L_\infty (M\dual(-1), -s)$ are multiplicative for $M \in \DMBeic(\Z)$, the latter up to sign \cite[Prop. III.1.2.8]{FPR}. 
By \refre{truly.geometric}\refit{truly.geometric.generators}, it is enough to show the claim for $M = i_* N$, $N \in \DMBeic(\Fp)$ and $M = E  := \eta_{!*} \eta^* \h^{-b}(X, -m)$, where $X / \Z$ is regular, flat and projective. For $M = E$, we have $L(M, s) = L_\Q(\h^{-b-1}(X_\eta), s)^{-1}$. This and the formula \refeq{Edual} for $M\dual(-1)$ in this case shows that the conjectural (see \ref{conj_genpropLseries}) properties of $L_\Q(\h^{-b-1}(X_\eta), s)$ implies the same properties for $L(M, s)$.
The $L$-series of $M=i_* N$ is a rational function in $p^{-s}$ (a priori with complex coefficients), which immediately yields the convergence for $\Re(s) \gg 0$ and the meromorphicity. 
\ifpreprint
We now show the functional equation for motives $i_* N$.\footnote{This must be well-known, but I failed to find a reference for it.} The $\ell$-adic realization commutes with duals: $(N\dual)_\ell = (N_\ell)\dual$. Therefore, we have to see
\eqn \mylabel{eqn_Lduality}
L(V, s) = a b^s L(V\dual, -s)
\xeqn
with some nonzero numbers $a, b$, for any finite-dimensional continuous complex representation $V$ of $\Gal(\Fp)$. Here $L(?, s) := \det(\Id-\Fr^{-1} p^{-s}|?)$. We may replace $V$ by $\det V$ without changing either side of the \refeq{Lduality}, so we may assume $\dim V = 1$. Then $\Fr^{-1}$ acts on $V$ ($V\dual$) by multiplication with some $f \in \CC^\x$ ($1/f$, respectively). Hence we can take $a:=-f$ and $b := 1/p$ in \refeq{Lduality}.
\fi
\iffinal
Noting that $(i_* N)\dual \twi{-1} = i_* (N\dual)$, the functional equation also holds unconditionally, as is well-known. 
\fi
\xpf

\bem \mylabel{bem_rationalfunceq}
Under \refco{numrat} the constant $a$ above is rational for $M = i_* N$, where $i: \Spec \Fp \r \Spec \Z$. To see this, we may assume by triangulatedness that $N$ is a pure motive with respect to numerical or homological equivalence, so that its $L$-function is a rational function in $p^{-s}$ with rational coefficients (see the reference in the proof of \refth{Fp}). 
\xbem

\section{Is the Arakelov intersection pairing perfect?}
\conj \mylabel{conj_motcomp}
For any geometric motive $M$ over $\Z$, $M \in \DMBeic(\Z)$ (see \refsect{motives} for the notation), the Arakelov intersection pairing between motivic homology and Arakelov motivic cohomology (\refde{A.i.p.})
\eqn \mylabel{eqn_A.i.p.}
\pi_M : \H_{-2}(M, -1)_\RR \x \Hhat^0_\RR(M) \r \RR
\xeqn
is a perfect pairing of finite-dimensional $\RR$-vector spaces. 
\xconj 

\rema \mylabel{rema_bem1}
\begin{enumerate}[(i)]
 \item 
The shape of \refeq{A.i.p.} is similar to the situation of \'etale constructible sheaves over $\SpecZ$: thinking of $M \in \DMBeic(\Z)$ as being analogous to a complex of constructible sheaves $\mathcal F$ over $\Z$, the groups $\HD^*(M)$ correspond (in spirit) to the Tate cohomology groups $\H^*_{\text{Tate}}(\RR, \mathcal F|_\RR)$ at the archimedean place. Given that, $\Hhat^i(M)$ parallels $\Hc^i(\mathcal F):= \H^i \RGc (\Z, \mathcal F)$, that is to say, cohomology with compact support, which is defined via $\RGc := \cone [-1] \left (\RG(\Z, \mathcal F) \r \RG_\text{Tate}(\RR, \mathcal F|_\RR) \right)$, much the same way as \refeq{onehat}, \refeq{Hhat}. Finally, the Arakelov intersection pairing corresponds to the perfect pairing of \emph{Artin-Verdier} duality, see e.g. \cite[Ch. II.3]{Milne:Arithmetic}
$$\Hc^{i}(\Z, \mathcal F) \x \Ext^{3-i}_\Z(\mathcal F, \Gm) \r \Hc^3 (\Z, \Gm).$$
A higher-dimensional extension was conjectured by Milne \cite[Conjecture II.7.17]{Milne:Arithmetic} and proven by Geisser \cite{Geisser:Duality}. 
\item \mylabel{item_motcomp.duality}
For any fixed $M \in \DMBeic(\Z)$, \refco{motcomp} for all $M[k]$ ($k \in \Z$) is equivalent to the one for $M\dual\twi{-1}[k]$. This follows from \refre{Hhat.basic}\refit{duality}, \refit{duality2} and the five lemma. 
\item
Gillet and Soul\'e conjecture that the intersection product 
\eqn \mylabel{eqn_pairingGS}
\CHtilde[m](X)_\RR \x \CHtilde[d-m](X)_\RR \r \RR
\xeqn 
is non-degenerate for any regular scheme $X$ that is projective and flat over $\Z$ of constant dimension $d$ \cite[Conjecture 1]{GilletSoule:Analogs}. 
By \refex{AipX}, at least for $X$ smooth, this pairing is compatible with the Arakelov intersection pairing $\pi_{\M(X)\twi{m-d}}$, i.e., there is a commutative diagram of pairings,
$$\begin{array}{cccccccccc}
0 & \r & \Hhat^0(M) = \CHhat^m(X)_\RR & \r & \CHtilde[m](X)_\RR & \stackrel \omega \r & \im \omega & \r & 0 \\
 &  & \x &  & \x & & \x & &  \\
0 & \leftarrow & \H_{-2}(M, -1) = \CH^{d-m}(X)_\RR & \leftarrow & \CHtilde[d-m](X)_\RR & \leftarrow & \im a & \leftarrow & 0 \\
 &  & \downarrow &  & \downarrow & & \downarrow & &  \\
 &  & \RR &  & \RR & & \RR & &  \\
  \end{array}
$$
where $\omega: \CHtilde[m](X) \r A^{m,m}(X)$ and $a: A^{d-m-1,d-m-1}(X) / (\im \partial + \im \ol \partial) \r \CHtilde[d-m](X)$ are defined in \cite[Section 3.3.4]{GilletSoule:Arithmetic}. 
I don't know whether the pairing on the right is a non-degenerate pairing, so the relation of Gillet-Soul\'e's conjecture and \ref{conj_L} is unclear. 
Note that $\im \omega$ and $\im a$ are infinite-dimensional $\RR$-vector spaces.
\end{enumerate}
\xrema


Next, we show that \refco{motcomp} recovers all the axioms on mixed motives over $\Fp$ we were willing to assume. Previously, it was known that Tate's conjecture about the pole order of $\zeta$-functions over finite fields and \refco{numrat} together imply the Beilinson-Parshin conjecture \cite[Thm. 1.2.]{Geisser:Tate}, and that the Beilinson-Parshin conjecture is equivalent to Bondarko's weight functor $\DMgmeff(\Fp) \r \K^\bound(\M_\rat^\eff)$ between the triangulated category of effective motives with the bounded homotopy category of effective Chow motives (with rational coefficients) being an equivalence of categories \cite[Section 8.3.2]{Bondarko:Differential}.

\theo \mylabel{theo_reductionperfFp}
\refco{motcomp} for motives of the form $M = i_* N$ ($N$ any geometric motive over $\Fp$, $i: \SpecFp \r \SpecZ$) is equivalent to the conjunction of \refco{numrat} and the \emph{Beilinson-Parshin conjecture} stating
\eqn \mylabel{eqn_BeilinsonParshin}
K_r(X)_\Q = 0
\xeqn
for any smooth projective variety $X$ over $\Fp$ and all $r > 0$. 

Under the axioms concerning the existence and cohomological dimension of mixed motives over $\Fp$ and the weight formalism (see \refax{mixed}), \refco{motcomp} for all motives $i_* N$ is equivalent to \refco{numrat}.
\xtheo

\pf
Using the axioms about mixed motives, we first show that \refco{numrat} implies the perfectness. By construction, cf. \refeq{Hhat}, $\Hhat^*_\RR(i_* N) = \H^*(N)_\RR$. By \refaxcohomdim, the cohomological dimension of $\DMBeic(\Fp)$ is zero, so that $\H^j(N) = \H^0(\pH^j N)$ and similarly for $N\dual$. By the same axiom, only finitely many $j$ yield a non-zero term. Therefore, we may replace $N$ by $\pH^j N$ and assume that $N$ is a mixed motive. Using the weight filtration we reduce to the case where $N$ is a pure motive. Under \refco{numrat}, all adequate equivalence relations agree, so we may regard $N$ as a Chow motive or as a pure motive with respect to numerical equivalence. By the semi-simplicity of pure numerical motives there is a decomposition $N = \one^r \oplus R$, where $R$ satisfies $\H^0_{\DMBeic(\Fp)}(R\dual) = \H^0_{\DMBeic(\Fp)}(R) = 0$. By functoriality of the pairing we get a commutative diagram
$$\begin{array}{rccccl}
\H^0(N)_\RR & \x & \H^0(N\dual)_\RR & \lr & \RR \\
\downarrow \cong & & \uparrow \cong & & \downarrow = \\
\H^0(\one^r)_\RR  & \x & \H^0(\one^r)_\RR & \lr & \RR \end{array}
$$
The lower line is a perfect pairing, since the one for $\one_{\Fp}$ is by \refex{Fp}. 

We now show the second statement. Let $X$ be a smooth equidimensional projective variety over $\Fq$. We regard it as a $\Z$-scheme. 
By \refex{Fp}, the Arakelov intersection pairing
$$\Khat_{2m-k}(X)^{(m)} \x K_{k-2m}(X)^{(\dim X-m)}_\Q= K_{2m-k}(X)^{(m)}_\RR \x  K_{k-2m}(X)^{(\dim X-m)}_\RR \r \RR$$
is the usual multiplication on Adams eigenspaces in $K$-theory, followed by the multiplication with $\log p$ (which is irrelevant for the question of the perfectness). For $2m - k > 0$ the second factors vanishes, hence the perfectness is equivalent to \refeq{BeilinsonParshin}. For $2m = k$ is perfectness is equivalent, by definition, to the agreement of numerical and rational equivalence (up to torsion). This shows one implication of the second statement. By resolution of singularities, the category $\DMBeic(\Fp)$ is generated as a thick category by motives $\M(X)(m)$ as above. Since the perfectness only has to be checked on such generators, we are done with the converse implication as well.
\xpf

The following corollary was pointed out to me by Bruno Kahn.
\coro The perfectness of $\pi_M$ for all motives $M = i_* N$ implies a canonical equivalence $\DMBeic(\Fp) = \D^\bound(\PureMot_\rat(\Fp))$, which in turn implies among other things the independence of $L$-functions of $\ell$. 
\xcoro
\pf
That description of $\DMBeic(\Fp)$ is a consequence of $\sim_\num = \sim_\rat$ and the Beilinson-Parshin conjecture \cite[proof of Theorem 56]{Kahn:AlgebraicKTheory}.
\xpf

We now give some interesting consequences of \refco{motcomp} for motives which are truly motives over $\Z$, i.e., not coming from a motive over $\Fp$. It would be interesting to know whether other axioms on mixed motives over $\Q$, such as the agreement of homological and numerical equivalence on smooth projective varieties $X_\eta / \Q$ can be derived from \refco{motcomp}.

\theo 
\mylabel{theo_Beilinson.Soule}
As in \refex{AipX}, consider the motive $M = \M(X)\twi {m-d}[p-2m]$, $X / \Z$ regular, flat, projective and of equidimension $d$. 
Then \refco{motcomp} for $M$ is equivalent to the \emph{Beilinson-Soul\'e vanishing conjecture}
$$K_{2m-p}(X)^{(m)}_\Q = 0 \ \ (\text{for }p < 0\text{ and for }p = 0,\ m > 0).$$
\xtheo
\pf
The group $\Hhat^0(M)$ appears in the long exact sequence
$$\dots \r \HD^{-1}(M) = \HD^{2d-p-1}(X, d-m) \r \Hhat^0(M) \r \H^0(M) = \underbrace{K_{p-2m}(X)^{(d-m)}}_{=0} \r \dots$$
where the right hand vanishing is because $p-2m < 0$ for $p < 0$ and $p = 0$, $m > 0$.
The left hand vector space is dual to $\HD^{p}(X, m)$ by \refeq{dualityweak} (note that $d = \dim X_\CC + 1$). It vanishes for $p < 0$ for trivial reasons. For $p = 0$, the short exact sequence \refeq{longDeligne2} gives $\HD^0(X, m) = 0$ for $ m > 0$. Indeed, the Hodge structure on $\HdR^0(X)$ only lies in the $(0,0)$-part of the Hodge diamond, i.e., $F^m = 0$ for $m > 0$. Hence the injectivity of $\HBetti^0(X, \RR(m)) \r \HBetti^0(X, \CC) \stackrel \cong \r \HdR^0(X_\CC)$ gives the claim.
Therefore \refco{motcomp} for $M$ is equivalent to $\H_{-2}(M, -1) = \H_{2m-p-2}(\M(X)\twi {m-d}, -1)  = K_{2m-p}(X)^{(m)}_\RR = 0$.
\xpf

\exam \mylabel{exam_specialHasseWeilI}
Using the notation of \refth{Beilinson.Soule}, the group $\H_{-2}(M, -1)$ vanishes for $2m-p < 0$. Therefore, \ref{conj_motcomp} asserts  
that the Chern class map 
\eqn \mylabel{eqn_Chernclassisoinj} 
\H^{0}(M)_\RR = K_{p-2m}(X)^{(d-m)}_\RR \r \HD^{0}(M) = \HD^{2d-p}(X, d-m)
\xeqn
is injective for $p-2m > 0$ and an isomorphism for $p - 2m > 1$. In particular, the non-torsion part of higher $K$-theory of $X$ is finitely generated---a weakening of \refco{Bass}. 
\xexam

\prop \mylabel{prop_cohomdim2}
Assuming the existence of motivic $t$-structure on $\DMBeic(\Z)$ such that Betti and de Rham realization are exact (see \refax{mixed}\refit{exactness} and \refeq{exactnessBetti}), the perfectness of the Arakelov intersection pairing for all motives $M \in \DMBeic(\Z)$ implies that the cohomological dimension of mixed motives over $\Z$ is two, i.e., $\Hom(\one[1], M[n]) = 0$ for any $n > 2$ and $M \in \MM(\Z)$.  
\xprop
\pf
Let $M$ be a mixed motive. The group $\H_{1-n}(M, -1)_\RR = \Hom(\one[1-n](-1), M)_\RR$ is zero for $n < 0$: in this case $\one[1-n]$ lies in degree $n < 0$ (with respect to the motivic $t$-structure). On the other hand this group is dual, via $\pi_{M[n-3]}$, to $\Hhat_\RR^{3-n}(M)$. In \refeq{Hhat}, this group lies between $\H^{3-n}_\RR(M)$ which vanishes for $n > 2$ for the same reason and the Deligne cohomology group $\HD^{3-n}(M) = \Hom(M, \HD[3-n])$ which in turn vanishes by exactness of Betti and de Rham realization, except for $n=1, 2$, as in \refeq{Delignemixed}. Consequently, $\H_{1-n}(M)=0$ except for $n=0$, $1$, $2$.
\xpf

\lemm \mylabel{lemm_finitelymany}
Under \refco{motcomp}, $\H^i(M)$ is nonzero only for finitely many $i \in \Z$.
\xlemm
This is a consequence of the spectral sequence $\H^a (\pH^b (M)) \Rightarrow \H^{a+b} (M)$, the boundedness of the motivic $t$-structure and of the cohomological dimension \refaxcohomdim. It also follows from the perfectness of the Arakelov intersection pairing (not using the axioms on mixed motives):

\pf
It suffices to check the claim for $M = \M(X)(m)$, where $X$ is as in \refex{specialHasseWeilI} and $m \in \Z$, since these objects generate $\DMBeic(\Z)$ as a thick category by resolution of singularities. Now, the claim follows as in \refpr{cohomdim2} using the vanishing $K_k(X)$ for $k < 0$ and the vanishing of almost all Deligne cohomology groups of $X$. 
\xpf

\section{Are special $L$-values given by the Arakelov intersection pairing?} \mylabel{sect_specialL}

Throughout this section, let $M$ be any geometric motive over $\Z$. In this chapter, wherever ranks of motivic cohomology groups are involved, we assume that the Bass conjecture holds up to torsion: 

\conj \mylabel{conj_Bass}
For any regular scheme $X / \Z$, $\dim_\Q K_i(X)_\Q < \infty$. 
\xconj

We need the following consequence (by resolution of singularities): motivic cohomology of all geometric motives over $\Z$ is finitely generated. 

By \refaxcohomdim\ (see also \refle{finitelymany}) and \refeq{Hhat} only finitely many $\H^i(M)$ and $\Hhat^i(M)$ are nonzero as $i \in \Z$ varies. Thus, the Euler characteristic 
\eqn
\mylabel{eqn_Eulerchar}
\chi(M) := \sum_i (-1)^i \dim \H^i(M)
\xeqn 
and similarly $\widehat \chi(M)$, $\chi_\Deligne(M)$ are well-defined. We write $\det H^* := \t_{i \in \Z} \det^{(-1)^i} H^i$ for any bounded family $H^i$ of finite-dimensional vector spaces, such as $\H^i(M)$ etc.\
The determinant of Arakelov motivic cohomology groups carries a $\Q$-structure by the isomorphism induced by \refeq{DeligneQ} and \refeq{Hhat},
$$\det \Hhat_\RR^*(M) = \left ( \det \H^*(M) \t \det^{-1} \H_{\Betti, \Q}^*(M) \t \det \HdR^{/F, \alg, *}(M) \right ) \t_\Q \RR.$$
 
\conj \mylabel{conj_L}
The order of the $L$-function of $M$ (\refde{LfunctionoverOF}) is given by
$$\ord_{s=0} L(M, s) = - \chi(M\dual(-1)).$$
As usual, negative orders mean a pole, positive ones a zero of the $L$-function. Moreover, assuming the perfectness of the Arakelov intersection pairings $\pi_{M[k]}$ (\refde{A.i.p.}) for all $k \in \Z$ asserted by \refco{motcomp}, the special $L$-value is given by 
$$L^*(M, 0) \equiv 1 / \Pi_{M} \mod {\Q^\x}.$$
Here $\Pi_M$ means the following: the perfectness of the Arakelov intersection pairing yields a map  
$$\det \H_{-2+*}(M, -1)_\RR  \t \det \Hhat_\RR^*(M) \r \RR.$$
The $\Q$-structure on the left maps to a real number denoted $\Pi_M$. Note that $\Pi_M$ is well-defined up to multiplication by a non-zero rational number.
\xconj

\nota \mylabel{nota_E}
For a projective flat scheme $X / \Z$ with smooth generic fiber $X_\eta / \Q$, we write $E := \eta_{!*} \eta^* \h^{-b}(X, -m) \in \MM(\Z)$ and $M_\eta = \eta^*[-1] E = \h^{-b-1}(X_\eta, -m) \in \MM(\Q)$. The definition of $E$ is recalled in \refsect{motives}. In particular, whenever $E$ is considered, we need to assume the axioms on mixed motives mentioned in \refsect{motives}. The motive $E$ only depends on $X_\eta$, not on $X$. It is pure of weight $w :=  \wt(E) = 2m - b$. Putting $d := \dim X$ and $d_\eta = \dim X_\eta$, the dual 
\eqn \mylabel{eqn_Edual}
E\dual = (\eta_{!*} \eta^* \h^{-2d+4+b}(X, 1-d+m))[-2]
\xeqn
is pure of weight $-w$, while $M_\eta$ is pure of weight $w-1$. 
\xnota

Under \refco{motcomp}, the pole order conjecture is equivalent to 
$$\ord_{s=0} L(M, s)  = - \widehat \chi (M).$$ 
We expound some structural properties of the conjecture. In order to state the compatibility with the functional equation, we shall need the following conjecture due to Deligne. It implies the compatibility of the $L$-values conjecture for critical pure motives $M_\eta$ over $\Q$ (i.e., motives such that $\Hw^i(M_\eta) = 0$, $i =0, 1$) with the functional equation \cite[Theorem 5.6]{Deligne:Valeurs}.

\conj \cite[Conjecture 6.6]{Deligne:Valeurs} \mylabel{conj_rankone}
Let $M$ be a pure motive over $\Q$ with respect to homological equivalence, i.e., a direct summand in $\PureMot_\hom(\Q)$ of $h(X_\eta, m)$ where $X_\eta / \Q$ is smooth projective. Assume that $M$ is of rank one, that is to say, its Betti realization (or, equivalently, de Rham or $\ell$-adic realization) is one-dimensional. Then $M$ is of the form $M(\epsilon)(n)$, where $n$ is an integer and $\epsilon: \Gal(\Q) \r \Q^\times$ is a finite character and $M(\epsilon)$ denotes the Dirichlet motive attached to the one-dimensional representation, $\epsilon$, of $\Gal(\Q)$ (\lc).
\xconj

\theo \mylabel{theo_structure}
\begin{enumerate}[(i)]
\item \mylabel{item_1}
\refco{L} is triangulated: given a distinguished triangle $M_1 \r M_2 \r M_3$ in $\DMBeic(\Z)$, the conjecture predicts
$$L^*(M_1, 0) L^*(M_3, 0) = L^* (M_2, 0)$$
and additively with the pole orders. In particular, the subcategory of $\DMBeic(\Z)$ of motives for which the conjecture holds is triangulated.
\item \mylabel{item_2}
Assume Deligne's \refco{rankone}, \refco{numrat} ($\sim_\rat = \sim_\num$), the functional equation for completed $L$-functions over $\Q$ (\refco{genpropLseries}) and \ref{conj_motcomp}. Then \refco{L} for any accessible motive $M$ (\refde{accessible}) is equivalent to the one for $M\dual\twi{-1}$. 
\end{enumerate}
\xtheo

Note that accessible motives generate $\DMBeic(\Z)$ as a thick category (\refre{truly.geometric}).

\pf \refit{1}: The pole order additivity is clear. The multiplicativity of the special values formula follows easily by considering the long exact sequences made of $\Hhat_\RR^*(M_i)$ and $\H_*(M,-1)_\RR$. By construction, the $\Q$-structure on Arakelov motivic cohomology is triangulated, i.e., there is a canonical isomorphism
$\det \Hhat_\RR^*(M_2) = \det \Hhat_\RR^*(M_1) \t \det \Hhat_\RR^*(M_3)$ of $\RR$-vector spaces, respecting the $\Q$-structure. 

\refit{2}: 
By \refre{truly.geometric}\refit{truly.geometric.generators}, it is enough to show the claim for all $M$ contained in the triangulated subcategory of $\DMBeic(\Z)$ generated by the image of $i_*: \DMBeic(\Fp) \r \DMBeic(\Z)$ for all primes $p$ and motives $E$ as in \refno{E}. 

We put $\ord := \ord_{s=0}$ and $\chi_w^a (M) := \sum_{i \in \Z} (-1)^i \dim \H^a_w(\H^i(\RG_\mathrm{H}(M))$ for $a = 0$, $1$, where $\RG_\mathrm{H}$ denotes the Hodge realization functor defined in \refeq{Hodgerealization}. \refco{L} for $M$, $\ord L(M, s) = - \chi(M\dual(-1))$, is equivalent to 
\eqnarr
\ord \Lambda(M) & \stackrel{\text{\refeq{ordLinfty}}} = & \ord L(M) - \chi_w^1 (M\dual(-1)) \\  
& = & - \chi(M\dual(-1)) - \chi^1_w(M\dual(-1)) \\
& \stackrel{\text{\ref{conj_motcomp}}}=  & - \widehat \chi(M) - \chi_w^1(M\dual(-1)) \\
& \stackrel{\text{\refeq{dualityweak}}} = & - \chi(M) + \chi_\Deligne(M) - \chi_w^0(M) \\
& = & - \chi(M) - \chi_w^1(M) 
\xeqnarr
Indeed, $\chi_\Deligne(M) = \chi_w^0(M) - \chi_w^1(M)$ (at least) for all $M$ as in the claim: for $M = E$, this follows from \refeq{specweak},  \refeq{Deligneweak}, and \refeq{Delignemixed}, while for $M = i_* N$, these terms are zero. 
By \refle{functionalequation},  the functional equation for mixed motives over $\Q$ implies the one for motives over $\Z$, so that $\ord \Lambda(M\dual(-1)) = \ord \Lambda(M)$. Again invoking the pole order calculation of $L_\infty$-functions we get $\ord L(M\dual(-1)) = - \chi (M)$, that is, the conjectural prediction of the pole order of $L(M\dual(-1))$. This settles the compatibility of the pole order prediction with the functional equation.

As for the special $L$-values, the claim is again triangulated. For motives $M = i_* N$,  where $i: \SpecFp \r \SpecZ$ and $N$ is any geometric motive over $\Fp$ we have $M\dual \twi{-1} = i_* N\dual$. The functional equation reads $L(i_* N, s) = a b^s L(i_* N\dual, -s)$, with $a$ and $b$ in $\Q^\x$ (\refbe{rationalfunceq}; this uses the agreement of numerical and homological equivalence, implied by \refco{numrat}.) On the other hand $\Pi_{i_* N} \equiv \Pi_{i_* N\dual} \text{ mod } \Q^\x$ by \refex{Fp}. 

To check the claim for $M = E$ as above, we may assume $X$ is of equidimension $d$. Each individual Deligne cohomology group $\HD^r(E)$ carries the $\Q$-structure mentioned in \refeq{longDeligne}. We can assume $w := \wt(E) = 2m-b \leq 2$, since otherwise we can replace $E$ by $E\dual(-1)[2]$. Let $M_\eta := \eta^* M[-1] = \h^{-b-1}(X_\eta, m)$, where $X_\eta$ is the generic fiber of $X$. For $w = 2$, the hard Lefschetz axiom implies an isomorphism $E \cong E\dual(-1)[2]$ (see \refeq{hardLefschetz}), so that there is nothing to show in that case. Let now $w \leq 1$. We write $L^*(-) := L^*(-,0)$. Deligne's conjecture \ref{conj_rankone} implies (see \lc)
$$\frac{L^*(E)}{L^*(E\dual(-1))} \stackrel{\text{\ref{satz_corrj}}}= \frac{L^*(M_\eta\dual(-1))}{L^*(M_\eta)}  \stackrel{\text{\ref{conj_rankone}}}\equiv \frac{a_1}{a_2} \text{ mod } \Q^\x$$
where $a_1$ is an element in the $\Q$-lattice of $\det \HD^0(E\dual(-1)) (= \RR)$ given by the $\Q$-structure on this Deligne cohomology group, and $a_2$ is an element in the $\Q$-lattice of $\det^{-1} \HD^{1}(E)$, regarded as an element of $\det \HD^0(E\dual(-1))$ using the isomorphism $\HD^1(M_\eta\dual(-1))\dual \r \HD^0(M_\eta)$, cf.\ \refeq{dualityweak}. In other words, the isomorphism $\det \HD^1(E) \r \det \HD^0(E\dual(-1))$ is multiplication by $a_1 / a_2$ with respect to the $\Q$-structures on both sides. 

For $r \neq 1$, the group $\HD^r(E)$ and its $\Q$-structure is trivial, since the corresponding Betti and (truncated) de Rham cohomology groups vanish. Therefore there is a canonical isomorphism 
\eqn \mylabel{eqn_QstructureHD}
\det \HD^*(E) \cong \det^{-1} \HD^1(E)
\xeqn 
(including the $\Q$-structure). 
Thus
\eqnarr
\det \Hhat_\RR^*(E) \t \det \H_*(E, -1)_\RR & =  & \det \H^*(E)_\RR \t \det^{-1} \HD^*(E) \t \det \H^*(E\dual, 1)_\RR \\
& \cong & \det \H^*(E)_\RR \t \det^{-1} \HD(E\dual(-1)) \t \det \H^*(E\dual, 1)_\RR \\
& = & \det \H_*(E\dual)_\RR \t \det \Hhat_\RR^*(E\dual,1) 
\xeqnarr
Both the left hand side and the last term on the right hand side map to $\RR$ via the Arakelov intersection pairings for $E$ and $E\dual\twi{-1}$, respectively. The two pairings are compatible with the isomorphism by the commutativity of \refeq{dualitycompatible}. By \refco{L} for $E$, the image of the $\Q$-structure on the left hand side is $L^*(E)^{-1}$, while the one from the right hand side is, by \ref{conj_L}, just $1 / L^*(E\dual(-1))$. Hence the two cases of the conjecture are equivalent.  
\xpf

In the remainder of this paper, we show how certain special cases of \ref{conj_L} are related to conjectures of Beilinson, Soul\'e, and Tate. In order to formulate our main result as succinctly as possible, we formulate the following
\conj \mylabel{conj_heightcomp}
For the motive $E$ defined in \refno{E} with $w := \wt(E) = 2$, the Arakelov intersection pairing $\pi_{E[-2]}: \H_0(E, -1) \x \Hhat^2(E) \r \RR$ agrees with Beilinson's height pairing \refeq{Beilinsonb=2m}.
\xconj 

\rema \mylabel{rema_height}
By \refth{summaryH1f} and \refeq{Delignemixed}, we know $\H_0(E,-1) = \CH^{d-m}(X_\eta)_\Q$ and $\Hhat^2(E) = \H^2(E) = \CH^{m}(X_\eta)_\Q$ (cf.\ the proof of \refpr{comparison}), so this conjecture only concerns the pairing itself. Moreover, \refeq{Beilinsonb=2m} is induced by the Gillet-Soul\'e intersection pairing 
$$\CH^{d-m}(X)_\Q \x \CHhat^{m}(X)_\Q \r \CHhat^d(X)_\Q \stackrel{f_*}\lr \CHhat^1(\Z) = \RR,$$
which in turn is induced by \refeq{pairingGS}. As mentioned in \refex{AipX}, this pairing agrees with the Arakelov intersection pairing for $\M(X)\twi{-m}$ at least up a $\Q$-linear combination of $\log p_i$, where $p_i$ are the primes such that the restriction of $X$ is smooth over $\Z[1/\prod p_i]$. It is worth mentioning that this comparison is an entirely formal consequence of the use of stable homotopy category. Its definition as the homotopy category of spectra of simplicial presheaves on \emph{smooth} schemes yields immediate comparison results such as \cite[Thm. 7.4]{Scholbach:ArakelovII} for smooth schemes, but not easily for other schemes. Therefore, it is a natural idea to overcome this hurdle by studying (Arakelov) motivic cohomology for log-smooth schemes. By de Jong's resolution of singularities, motives of all log-smooth $\Z$-schemes should generate a category of motives of logarithmic schemes over $\Z$. This would allow to bypass \refco{heightcomp}. I plan to return to this question in a 
subsequent paper.   
\xrema

The following two theorems summarize the remainder of this paper: under standard assumptions on motives and their $L$-functions, it shows that Beilinson's, Soul\'e's, and Tate's conjectures are essentially equivalent to the conceptual reformulation made possible by the use of the Arakelov intersection pairing.

\theo \mylabel{theo_summary.poles}
The following are equivalent:
\begin{enumerate}[(i)]
\item \mylabel{item_77.a}
The conjecture of Soul\'e (\ref{conj_SouleII}), restricted to regular, projective (but not necessarily flat) schemes. 
\item \mylabel{item_88.a}
The restriction of the pole order formula (\refco{L}) to the category of easily accessible motives (\refde{accessible}).
\end{enumerate}
\xtheo

\pf
This follows immediately from \refth{comparisonSoule} by \refth{structure}.
\xpf

By \refre{truly.geometric}, the thick closure of the category of easily accessible motives is the entire category $\DMBeic(\Z)$. Thus, the pole order formula of \refco{L} can be regarded as an extension of Soul\'e's conjecture to direct summands.

\theo \mylabel{theo_summary}
We assume the existence of mixed motives as formulated in \refax{mixed} and the agreement of Beilinson's height pairing with the Arakelov intersection pairing (\refco{heightcomp}). Moreover, in order to incorporate the compatibility of $L$-values with respect to the functional equation, we assume Deligne's conjecture \ref{conj_rankone} on rank one motives, and the functional equation for completed $L$-functions over $\Q$ (\refco{genpropLseries}). 
Finally, we assume that the pole order formula of \refco{L} holds for all motives in $\DMBeic(\Z)$.

Then, the following are equivalent:
\begin{enumerate}[(i)]
\item \mylabel{item_77}
The conjunction of the conjectures of Beilinson ($L$-values and $\sim_\num = \sim_\rat$, \ref{conj_Beilinson}, \ref{conj_numrat}), and Tate (\ref{conj_Tate}).
\item \mylabel{item_88}
The restriction of the conjunction of the perfectness of the Arakelov intersection pairings (\refco{motcomp}) and the special $L$-values formula (\refco{L}) to the subcategory $\DMtrgm(\Z) \subset \DMBeic(\Z)$ of accessible motives (\refde{accessible}).
\end{enumerate}
\xtheo

\pf
By \refre{truly.geometric}\refit{truly.geometric.generators}, $\DMtrgm(\Z)$ is contained in the triangulated category generated by motives $M = E$ as in \refno{E}, and motives of the form $M = i_* N$, $N \in \DMBeic(\Fp)$, $i: \SpecFp \r \SpecZ$. For the latter type of motives, \refco{motcomp} is equivalent to \refco{numrat} by \refth{reductionperfFp} and \ref{conj_L} is equivalent to the Tate conjecture by \refth{Fp}.

The subcategory of $\DMBeic(\Z)$ of motives $M$ for which all pairings $\pi_{M[k]}$ are perfect is triangulated since motivic and Arakelov motivic cohomology behave well under triangles. 
Moreover, \ref{conj_motcomp} for $M (\in \DMtrgm(\Z))$ is equivalent to \ref{conj_motcomp} for $M\dual\twi{-1}$ by \refre{bem1}\refit{motcomp.duality}.
In a similar vein, \refco{L} is stable under distinguished triangles, and \ref{conj_L} for $M$ is equivalent to \ref{conj_L} for $M\dual\twi{-1}$ (\refth{structure}).

To finish \refit{77} $\Rightarrow$ \refit{88}, using the calculation of $E\dual\twi{-1}$ in \refeq{Edual}, we therefore only need to consider $M = E$ with $w := \wt(E) = 2m - b \leq 2$. Beilinson's pole order conjecture for $M_\eta$, \ref{conj_Beilinson}\refit{Beilinson1}, is equivalent (see \refeq{BeilinsonEulerChar}) to 
\eqn \mylabel{eqn_Beilinson}
\ord_{s=0} L(E, s) = - \chi(E\dual(-1)) + \dim \H^1(E\dual(-1)).
\xeqn
By assumption, $L(E, s) = - \chi(E\dual\twi {-1}) = -\chi(E\dual(-1))$, so that we get $\H^1(E\dual(-1)) = 0$. Using this vanishing, part \refit{Beilinson2} of Beilinson's conjecture is equivalent to the perfectness of the intersection pairings $\pi_{E[k]}$, $k \in \Z$ (with $w = \wt(E) \leq 2$), by \refpr{comparison}. This shows that \ref{conj_numrat}, \ref{conj_SouleII}, and \ref{conj_Beilinson}\refit{Beilinson2} together imply \ref{conj_motcomp} for all $M \in  \DMtrgm(\Z)$. Then parts \refit{Beilinson1}, \refit{Beilinson3} of Beilinson's conjecture are equivalent to \ref{conj_L} for all motives of the form $E$ (of weight $\leq 2$), by \refth{comparisonBeilinson}.

The converse implication \refit{88} $\Rightarrow$ \refit{77} is shown using the same arguments.
\xpf

\rema 
\mylabel{rema_natural}
It is natural to ask for the equivalence of the following two statements:
\begin{enumerate}[(i)]
\item 
\mylabel{item_77.c}
The conjectures of Beilinson, Soul\'e, and Tate (\ref{conj_Beilinson}, \ref{conj_numrat}, \ref{conj_SouleII}, \ref{conj_Tate}).
\item
\mylabel{item_88.c}
The restriction of Conjectures \ref{conj_motcomp} and \ref{conj_L} to the category of accessible motives. 
\end{enumerate}
Under the assumptions of \ref{theo_summary}, \emph{except} for the pole order formula assumption, the above proof does show \refit{88.c} $\Rightarrow$ \refit{77.c}. 
The latter addditional assumption is only needed to prove the converse, and is actually only needed for motives of the form $M = E$ as above.
Moreover, it holds unconditionally if $\M(X_\eta)$ is an Artin-Tate motive (\refth{H1Tate}).
The vanishing $\H^1(E\dual(-1))=0$ also follows from the Soul\'e+Tate conjecture if one can show $E \in \DMtrgm(\Z)$, which in its turn would follow if the motivic $t$-structure on $\DMBeic(\Z)$ restricts to a $t$-structure on $\DMtrgm(\Z)$. In this case, the proof of \cite[Prop. 5.6]{Scholbach:fcoho} referred to in \refre{truly.geometric}\refit{truly.geometric.generators} could be adapted to $\DMtrgm(\Z)$.
\xrema

\subsection{Relation to a conjecture of Soul\'e} \mylabel{sect_relationSoule}

\conj \mylabel{conj_SouleII} (Soul\'e, \cite[Conjecture 2.2.]{Soule:KTheorie})
Let $Y / \Z$ be quasiprojective. Let $m \in \Z$ be arbitrary. Then
\eqn \mylabel{eqn_Soule}
\ord_{s=m} \zeta(Y, s) = \sum_{i \geq 0} (-1)^{i+1} \dim_\Q K'_i(Y)_{(m),\Q}
\xeqn
We refer to \lc\, for the definition of the Adams eigenspace $K'_i(Y)_{(m),\Q}$. For $Y$ regular, it agrees with $K_i(Y)^{(\dim Y-m)}_\Q$.
\xconj

Soul\'e's conjecture extends a previous conjecture of Tate \cite[p. 105]{Tate:Algebraic}. A formally similar conjecture was also expressed by Lichtenbaum \cite{Lichtenbaum:Values}. The right hand side of \refeq{Soule} makes sense under the Bass conjecture \ref{conj_Bass} and the vanishing of almost all $K'$-groups, which in turn is a consequence of \refaxcohomdim. See also \refle{finitelymany}. As the thick closure of $\DMtrgm(\Z)$ is all of $\DMBeic(\Z)$, the following statement can be paraphrased by saying that Soul\'e's conjecture is essentially equivalent to the pole order part of \refco{L}. This proof does not make use of mixed motives. 

\theo \mylabel{theo_comparisonSoule}
\refco{SouleII} for $Y$ and $m$ is equivalent to the pole order prediction of \refco{L} for  $M = \Mc(Y)(-m)$. 
\xtheo
\pf
\refsa{HasseWeil} says $\zeta(Y, s+m) = L(\Mc(Y)(-m), s)$. The statement for $Y$ is implied by the conjunction of the one for some open subscheme $U$ of $Y$ and $Z := Y \backslash U$, since Adams eigenspaces in $K'$-theory have a localization sequence \cite[1.3.]{Soule:KTheorie}, and motives with compact support behave well: $\Mc(Z) \r \Mc(Y) \r \Mc(U)$ is a distinguished triangle. In particular we may assume that $Y$ is integral. Thus, there is an open affine subscheme $U$ of $Y$ that is either smooth over $\Z$ or over some $\Fp$: if $Y / \Z$ is flat, one can take an open neighborhood of a smooth point of the generic fiber of $Y$, otherwise $Y$ lies over some $\SpecFp$ and one can take a neighborhood of a smooth point of 
\ifpreprint
$Y$ \cite[Prop. 17.15.12.]{EGA4}. 
\fi
\iffinal
$Y$.
\fi
Let $f : Y \r \Z$ be the projection. By Noetherian induction, we may replace $Y$ by $U$ and hence assume $Y$ is regular and affine of dimension $d$, so that $(\Mc(Y)(-m))\dual\twi{-1} = f_! f^* \one(m)\twi {-1} = f_! f^! \one(m) \twi{-d}$ by purity. Hence
\eqnarr
\chi ((\Mc(Y)(-m))\dual\twi{-1}) & = & \chi(\M(Y)(m-d)[-2d]) \\
& = & \sum_{i \in \Z} (-1)^i \dim \H^{i+2d}(Y, m-d) \\
& = & \sum_{i \in \Z} (-1)^i \dim K_i(Y)^{(d-m)}_\Q. 
\xeqnarr
\xpf

\exam \mylabel{bsp_specialHasseWeil}
We continue Examples \ref{exam_AipX} and \ref{exam_specialHasseWeilI} and look at the special values of the $\zeta$-function of $X$: by \refsa{HasseWeil} we have $L(M, s) = \zeta(X, s+d-m)$. The Arakelov intersection pairing $\pi_{M[i]}$ concerns the following groups
$$\begin{array}{cccc}
0 & \x & 0 & i \leq -1 \\
K_{0}(X)^{(m)}_\RR & \x & \Hhat^0_\RR(M) & i = 0 \\
K_{1}(X)^{(m)}_\RR & \x & \coker K_0(X)^{(d-m)}_\RR \r \HD^{2(d-m)}(X, d-m)  & i = 1 \\
K_{i}(X)^{(m)}_\RR & \x & \HD^{2(d-m)+i-1}(X, d-m) & i > 1.
\end{array}$$
The pairing for $i \geq 1$ is given by the Chern class $K_{i}(X)^{(m)}_\RR \r \HD^{2m-i}(X, m)$ together with the cup product on Deligne cohomology, followed by the push-forward $f_*: \HD^{2d-1}(X, d) \r \HD^1(U, 1) = \RR$. I expect that the group $\Hhat^0_\RR(M)$ is isomorphic to $\CH^{d-m}(X)_\RR$ and that the pairing $\pi_M$ is the natural pairing of (arithmetic) Chow groups (cf.\ \refre{height}). We do know that these two pairings agree up to a $\Q$-linear combination of $\log \prod p_i$, where $p_i$ are the primes such that the restriction of $X$ is smooth over $\Z[1/\prod p_i]$.  

These pairings assemble to a map 
$$\bigotimes \pi_{M[i]} : \bigotimes_i \det^{(-1)^i} (\H_{-2-i}(M,-1)_\RR \t \Hhat_\RR^{i}(M)) \stackrel{\cong}\r \RR.$$ 
(Even though the groups $\Hhat^i(M)$ vanish for $i < 0$, the determinant carries a non-trivial information related to these groups, namely the determinants of the Chern class map, see \refeq{Chernclassisoinj}.) \refco{L} asserts that---modulo $\Q^\x$---$L^*(M, 0)$ is the reciprocal of the image of $1$ in $\RR$ via the $\Q$-structure map of the left hand term. The class number formula has been interpreted in terms similar to the one above, see \cite[III.4.3]{Soule:Arakelov}.
\xexam

\subsection{Relation to Beilinson's conjecture} \mylabel{sect_relationBeilinson}

In this section, we use the notation of \ref{nota_E}. The following is Beilinson's conjecture \cite{Beilinson:Higher,Beilinson:Notes}. Part \refit{Beilinson1} concerns the pole order of $L$-functions, part \refit{Beilinson2} is about the relation of Deligne cohomology and motivic cohomology, and \refit{Beilinson3} expresses the special $L$-value up to $\Q^\x$ in terms of determinants of the isomorphisms asserted in \refit{Beilinson2}. The pole order conjecture in case $w=3$ is due to Tate \cite{Tate:Algebraic}.
\conj \mylabel{conj_Beilinson} Let $X_\eta / \Q$ be smooth and projective.
\begin{enumerate}[(A)]
\item \mylabel{item_Beilinson1}
$$\ord_{s=m} L_\Q(\h^{-b-1}(X_\eta), s) = \ord_{s=0} L_\Q(M_\eta, s) = \left \{ \begin{array}{ll}
0 & w \geq 4 \\
-\dim \CH^{n}(X_\eta)_\Q / \hom & w = 3 \\
\dim \CH^{n} (X_\eta)_{\Q,\hom} & w = 2 \\
\dim \H^{b+2} (X_\eta, n)_\Z & w \leq 1
\end{array} \right.$$
Here $n := b+2-m = m+2-w$, and the groups at the right have been defined in \refsect{motives}.
\item \mylabel{item_Beilinson2}
For $w=2$, the \Def{height pairing} 
\eqn  \mylabel{eqn_Beilinsonb=2m}
\CH^m(X_\eta)_{\Q,\hom} \t \CH^{d - m}(X_\eta)_{\Q,\hom} \r \RR
\xeqn
is perfect.  

For $w = 1$, the map 
\eqn \mylabel{eqn_Beilinsoni=2m}
r_\infty: (\CH^m(X_\eta)_\Q/\hom \oplus \H^{2m+1}(X_\eta, n)_\Z) \t_\Q \RR \r \HD^{2m+1} (X, n).
\xeqn
obtained by the composition
$$\CH^m(X_\eta)_\Q/\hom \t \RR \r \H^{2m}_\dR(X_\RR) \r \HD^{2m+1}(X, n)$$
(see \refeq{longDeligne} for the right hand map) and
the realization map, is an isomorphism. 

For $w \leq 0$, the statement is the same, except that \refeq{Beilinsoni=2m} gets replaced by 
\eqn \mylabel{eqn_Beilinsoni>2m}
r_\infty: \H^{b+2}(X_\eta, n)_\Z \t_\Q \RR \r \HD^{b+2} (X_\eta, n).
\xeqn
\item \mylabel{item_Beilinson3}
The special $L$-value $L^*(M_\eta, 0)$  is conjecturally given up to a nonzero rational multiple by the following:

For $w=2$, by the determinant of the height pairing \refeq{Beilinsonb=2m} multiplied with the period of $M_\eta$, that is to say, the determinant of the isomorphism 
$$\alpha_{M_\eta}: \H^{2m-1}(X_\eta(\CC), \RR(m))^{(-1)^m} \r \H^{2m-1}_\dR((X_\eta)_\RR) / F^m$$
with respect to the usual $\Q$-structures on both sides (compare \refeq{longDeligne}).

For $w=1$, the $L$-value is given, mod $\Q^\x$, by $d_\infty(1)$, where
$$d_\infty := \det r_\infty: \det (\H^b(X_\eta, m)_\Z \oplus \CH^m(X_\eta) / \hom)_\RR = \RR \r \det \HD^b (X, n) = \RR,$$
the left hand term is endowed with the obvious $\Q$-structure, the right one gets the one stemming from the identification of $\HD^b (X_\eta, n) = \Hw^1(\H^{b-1}(X, \Q(n)))$ with the dual of $\Hw^0(\H^{b-1}(X, \Q(n))\dual(1))$. 

For $w\leq 0$, the statement is the same, except that the term $\CH^m(X_\eta) / \hom$ is omitted.
\end{enumerate}
\xconj

This concludes the statement of Beilinson's conjecture. It predicts $L$-values of motives $\h^{-b-1}(X_\eta, m)$ with $w = 2m-b \leq 2$, up to a nonzero rational factor. The remaining weights are adressed by the functional equation (\refco{genpropLseries}).

We compare Beilinson's conjecture with Conjecture \ref{conj_motcomp} and \ref{conj_L} applied to the generic intermediate extension $E := \eta_{!*} \eta^* \h^{-b}(X, -m)$, where $X$ is any projective model of $X_\eta$ (see \refno{E}). 

Recall from \cite[5.4.2.1]{Andre:Motifs} that the agreement of homological and numerical equivalence (which is part of \refax{mixed}) implies the hard Lefschetz isomorphism:
\eqn \mylabel{eqn_hardLef}
\h^{-b-1}(X_\eta, m-b-2) \stackrel{\cong}{\lr} \h^{- 2d_\eta + b+1}(X_\eta, -d_\eta+m-1) = M_\eta \dual(-1).
\xeqn
For $b + 1 \leq d_\eta$ the map is given by the $(d_\eta-b-1)$-st power of cup product with a hyperplane section, with respect to some embedding $X_\eta \subset \P[N]_\Q$. The right hand term of \refeq{hardLef} is $M_\eta\dual(-1)$ by relative purity, applied to the smooth map $X_\eta / \Q$.  

\lemm \mylabel{lemm_Lefschetzprelim}
The hard Lefschetz isomorphism \refeq{hardLef} yields an isomorphism
\eqn \mylabel{eqn_hardLefschetz}
E\dual(-1)[2] = E(m-n) = E(w-2).
\xeqn
It induces isomorphisms of motivic and Deligne cohomology groups (respecting the $\Q$-structure of the latter):
\eqnarr
\CH^m(X_\eta)_\Q/\hom & \cong &  \CH^{d - m - 1}(X_\eta)_\Q/\hom,\\
\CH^m(X_\eta)_{\Q,\hom} & \cong & \CH^{d - m}(X_\eta)_{\Q,\hom} \text{ \cite[Conj. 5.3.(a)]{Beilinson:Height}},\\
\H^{b}(X_\eta, b-m)_\Z & \cong & \H^{2 d - b}(X_\eta, d-m)_\Z \text { for }w = 2m-b \leq 1 \ \  \\ 
\HD^{b}(X_\eta, b-m) & \cong & \HD^{2 d- b}(X_\eta, d-m) \text { for }w \leq 1. 
\xeqnarr
\xlemm

\pf
\refeq{hardLefschetz} is obtained from \refeq{hardLef} by applying $\eta_{!*}[1]$. Now apply \refth{summaryH1f} and the calculation of Deligne cohomology in \refeq{Delignemixed}.
%
\xpf

The following proposition compares the perfectness of certain Arakelov intersection pairings with the statements in part \refit{Beilinson2} in Beilinson's conjecture. 

\prop \mylabel{prop_comparison}
Let $E$ be as in \refno{E} with weight $w = \wt(E) \le 2$. If the weight of $E$ is $2$, we assume \refco{heightcomp}. The following are equivalent:
\begin{enumerate}[(i)]
\item The pairings $\pi_{E[i]}$ and $\pi_{E\dual\twi{-1}[i]}$ ($i \in \Z$) are perfect.
\item Part \refit{Beilinson2} of Beilinson's conjecture and $\H^3(E\dual(-1)[2]) = \H^1(E\dual, 1) = 0$ (only needed if $w \leq 1$).
\end{enumerate}
\xprop

\rema 
\mylabel{rema_vanish.ATM}
The group $\H^1(E\dual,1)$ vanishes unconditionally if $X_\eta$ is such that $M_\eta$ is a mixed Artin-Tate motive over $\Q$ (as opposed to a general mixed motive) by \refth{H1Tate}. 
Recall from \refth{summaryH1f} that $\H^3(E) = 0$ for $w := \wt(E) \leq 2$. 
\xrema

\pf
The proof combines the hard Lefschetz isomorphism (\refle{Lefschetzprelim}) and the calculation of motivic and Deligne cohomology of $E$ and $E\dual(-1)$ (\refth{summaryH1f}, \refeq{Delignemixed}).
  
The map $\H^{b+2}(X_\eta, n)_\Z \r \HD^{b+2}(X_\eta, n)$ featuring in \refeq{Beilinsoni=2m}, \refeq{Beilinsoni>2m} in the cases $w \leq 1$ of \refco{Beilinson} is the Chern class map $\H^2(E(m-n)) \r \HD^2(E(m-n))$. Via hard Lefschetz, this is the same as the Chern class map 
\eqn \mylabel{eqn_proofreal} 
\ch(E\dual(-1)): \H^0(E\dual(-1)) \r \HD^0(E\dual(-1)).
\xeqn  

Consider the case $w=1$. By Fontaine's reformulation \cite[9.5]{Fontaine:Valeurs}, the map \refeq{Beilinsoni=2m} being an isomorphism is equivalent to the existence of an exact sequence whose right hand map is the composition of the Poincar\'e duality isomorphism $\phi$ stemming from \refeq{dualityweak}, the hard Lefschetz isomorphism and the Chern class map. 
\small
$$
\xymatrix{
0 \ar[r] &
\CH^m(X_\eta)_\RR / \hom \ar[r]^{\ch} &
\HD^{2m}(X_\eta, m) \ar[r] \ar[d]_\phi^\cong &
\H^{2m+1}(X_\eta, m+1)_\Z \dual \t \RR \ar[r] &
0 \\
&
&
\HD^{2d_\eta-2m+1}(X_\eta, d_\eta+1-m)\dual \ar[r]^{\refeq{hardLef}}_\cong &
\HD^{2m+1}(X_\eta, m+1)\dual. \ar[u]_{\ch \dual}   &
}$$
\normalsize
In terms of motivic and Deligne cohomology groups, it reads
\eqn
\xymatrix{
0 \ar[r] &
\H^{1}(E)_\RR \ar[r]^{\ch^1(E)} &
\HD^{1}(E) \ar[r] \ar[d]_\phi^\cong &
\H^0(E\dual(-1))_\RR\dual \ar[r] & 
0 \\
&
&
\HD^0(E\dual(-1))\dual. \ar[ur]_{\ch^0(E\dual(-1))\dual}
} \mylabel{eqn_diag}
\xeqn
These groups also occur in the following exact sequences, whose terms are paired by the pairings indicated on top:
\eqn \mylabel{eqn_proof2}
\begin{array}{ccccccccc}
   \pi_{E[-1]}:   &    &\pi_{E\dual(-1)[-1]}: &&            &   & \pi_{E[-2]}: &  & \pi_{E\dual(-1)}: \\
\Hhat_\RR^{1}(E) & \r & \H^{1}(E)_\RR &  
\r & \HD^{1}(E) & \r & \Hhat_\RR^{2}(E) & \r & \H^{2}(E)_\RR  \\
\x & & \x &  & \x & & \x & & \x \\
\H^{1}(E\dual, 1)_\RR&\l& \Hhat_\RR^{1}(E\dual, 1) & \l & \HD^{0}(E\dual, 1)&
\l
&\H^{0}(E\dual, 1)_\RR & \l & \Hhat_\RR^{0}(E\dual, 1) 
\end{array}
\xeqn
We have $\H^2(E) = 0$, so the injectivity of $\ch^0(E\dual(-1))$ is equivalent to $\pi_{E\dual(-1)}$ being perfect. The identification of $\coker \ch^{1}(E)$ with $\H^0(E\dual,1)_\RR\dual$ of \refeq{diag} is equivalent to $\pi_{E[-2]}$ being perfect. The Chern class map $\H^1(E)_\RR = \CH^m(X_\eta)/\hom \r \HD^1(E) = \HD^{2m}(X_\eta, m) \subset \HBetti^{2m}(X_\eta, \RR(m))$ is injective by definition of homological equivalence. Hence $\Hhat^1_\RR(E) = 0$ so that $\H^1(E\dual,1)=0$ is equivalent to $\pi_{E[-1]}$ being perfect. By the five lemma, $\pi_{E\dual(-1)[-1]}$ is then perfect, too. All other Deligne, motivic, and hence Arakelov motivic cohomology groups of $E\dual(-1)$ and $E$, except for the ones displayed above, vanish. 

The case $w < 1$ is done similarly: in addition to the above, we have $\H^{1}(E) = 0$. Accordingly, \refeq{diag} reduces to an isomorphism $\HD^{1}(E) \stackrel{\cong}\r \H^0(E\dual(-1))_\RR\dual$. The details are omitted.

For $w = 2$, \emph{all} Deligne cohomology $\HD^*(E)$ and $\HD^*(E\dual, 1)$ vanish for weight reasons. Moreover $\H^a(E) = \H^{a-2}(E\dual(-1)) = 0$ for $a \neq 2$, so that $\pi_{E\dual(-1)[-1]}$ and $\pi_{E[-1]}$ are perfect.
The height pairing \refeq{Beilinsonb=2m} is just $\pi_{E[-2]}$ according to \refco{heightcomp}. Its perfectness is equivalent to the one of $\pi_{E\dual(-1)}$. 
\xpf

\theo \mylabel{theo_comparisonBeilinson}
We assume the perfectness of the Arakelov intersection pairing for motives of the form $M = E[n]$, with $E$ as in \refno{E} and $n \in \Z$.
We also assume \refco{heightcomp} if $E$ is of weight $2$. 
Then Beilinson's conjecture (parts \refit{Beilinson1}, \refit{Beilinson3}) for $M_\eta$ is equivalent to \refco{L} for $E$.
\xtheo 

\pf
By hard Lefschetz (\refle{Lefschetzprelim}) and calculation of motivic cohomology of $E$, \refth{summaryH1f}, part \refit{Beilinson1} of Beilinson's conjecture reads
\eqn \mylabel{eqn_BeilinsonEulerChar}
\ord_{s=0} L_\Q(M_\eta, s) \stackrel{\text{\ref{satz_corrj}}}= - \ord_{s=0} L_\Q(E, s) = \sum_{a \neq 1} (-1)^a \dim \H^a(E\dual(-1)).
\xeqn
In fact, $\H^a(E\dual(-1)) \stackrel{\text{\refeq{hardLefschetz}}}=\H^{a+2}(\eta_{!*} \eta^* \h^{-b}(X, -n))$. For example, in case $w = 2m-b \leq 1$, this equals $\H^{b+2}(X_\eta, n)_\Z$ for $a=0$ and vanishes for $a \neq 0, 1$. As was mentioned above, the perfectness of $\pi_{E[-1]}$ conjectured in \ref{conj_motcomp} implies $\H^1(E\dual(-1)) = 0$. (In case $w \geq 2$, we know this vanishing without invoking \ref{conj_motcomp}.) This settles the pole order part \refit{Beilinson1} of Beilinson's conjecture. 

For the special $L$-values, we revisit the proof of \refpr{comparison} and look at the involved $\Q$-structures. Again using hard Lefschetz, we replace the map $\H^{b+2}(X_\eta, n)_\Z \t \RR \r \HD^{b+2}(X_\eta, n)$ occurring in \refeq{Beilinsoni=2m}, \refeq{Beilinsoni>2m} by $\ch(E\dual(-1))$, see \refeq{proofreal}. The involved $\Q$-structures remain unchanged.  

We first treat the case $w = 1$. By \cite[9.5]{Fontaine:Valeurs}, \cite[Conj. III.4.4.3]{FPR}, Beilinson's conjecture is equivalent to saying that the $L$-value of $M_\eta$ is given by the reciprocal of the image (in $\RR$) of the $\Q$-structure on the right hand side:
\eqnarr
\RR & \cong  &  \det^{-1} \H^0(E\dual(-1))_\RR \t \det^{-1} \HD^{1} (E) \t \det \H^{1}(E)_\RR \\
& \stackrel{\text{\refeq{QstructureHD}}}= & \det^{-1} \H^*(E\dual(-1))_\RR \t \det \HD^*(E) \t \det^{-1} \H^*(E)_\RR.
\xeqnarr
Here the isomorphism stems from the exact sequence \refeq{diag} and the $\Q$-structure on $\HD^{1}(E)$ is the natural one defined in \refsect{absolute}. (This $\Q$-structure is distinct from the one on the isomorphic group $\HD^0(E\dual(-1))\dual$, as is apparent from the discussion of the functional equation in \refth{structure}.) Moreover, all groups $\H^*(E\dual(-1))$ except $\H^0$ and $\H^*(E)$ except $\H^1$ vanish. 
By construction of Arakelov motivic cohomology the above is isomorphic, including the $\Q$-structure, to
$$\det^{-1} \H^*(E\dual(-1))_\RR \t \det^{-1} \Hhat_\RR^* (E).$$
The above identification with $\RR$ agrees with the dual of the Arakelov intersection pairing for $E$, so that $L^*(E, 0) = L^*(M_\eta, 0)^{-1}$ is indeed the inverse of $\Pi_{E}$. This accomplishes the case $w = 1$.  

Again, the case $w \leq 0$ is similar but simpler, since in addition $\H^{1}(E)=0$. Correspondingly, only the determinant of the realization map $\ch(E\dual(-1))$ \refeq{proofreal}, as opposed to the one of \refeq{diag}, appears in Beilinson's conjecture. 

Finally, consider the special value at the central point, i.e., $w = 2$. In this case all groups $\HD^*(E)$ are trivial, but the $\Q$-structure on
$$\det^{-1} \HD^*(E) = \det \HD^*(M_\eta) = \det \underbrace{\H_\Betti^{b+1}(X_\eta, \RR(m))^{(-1)^m}}_{=:B} \t \det^{-1} \underbrace{\H_\dR^{b+1}(X_\eta \x \RR) / F^m)}_{=:dR}$$
is non-trivial since the period isomorphism $\alpha: B \r dR$ does not respect the natural $\Q$-structures. By linear algebra, $\det \alpha$ agrees (modulo $\Q^\x$) with the image (in $\RR$) of the $\Q$-lattice under the natural isomorphism induced by $\alpha$: $\det B \t \det^{-1} dR \stackrel \cong \r \RR.$ Except for $\H^2(E)$ and $\H^0(E\dual(-1))$, all motivic cohomology groups of $E$ and $E\dual(-1)$ vanish (\refth{summaryH1f}). The Arakelov intersection pairing $\pi_{E[-2]}$ agrees with the height pairing under \refco{heightcomp}. By \refeq{Hhat}, we have an isomorphism of $\RR$-vector spaces respecting the $\Q$-structure
$$\det \Hhat_\RR^*(E) = \det \H^*(E)_\RR \t \det^{-1} \HD^*(E),$$ 
so Beilinson's conjecture is indeed equivalent to saying that $L^*(E, 0) = L^*(M_\eta, 0)^{-1}$ is the reciprocal of the image of the $\Q$-lattice under $\det \Hhat_\RR^*(E) \t \det \H^*(E\dual(-1))_\RR \r \RR$.
\xpf

\subsection{Relation to the Tate conjecture over $\Fp$}\mylabel{sect_comparisonTate}


\conj  (Tate conjecture over finite fields \cite{Tate:Algebraic}) \mylabel{conj_Tate}
Let $X / \Fq$ be smooth and projective. Let $\ell$ be a prime such that $\ell \nmid q$. Any $\Gal(\Fq)$-invariant element of $\H^{2i}(X \x_{\Fq} \Fqq, \Ql(i))$ is a $\Ql$-linear combination of algebraic elements, i.e., elements in the image of the cycle class map $\CH^i(X) \r \H^{2i}(X \x_{\Fq} \Fqq, \Ql(i))$. 
\xconj

\theo \mylabel{theo_Fp}
In addition to the general assumptions on mixed motives over $\Fp$ (\refsect{motives}), we assume \refco{numrat}. Then the Tate conjecture \ref{conj_Tate} is equivalent to \refco{L} for motives $M = i_* N$, where $N$ is any geometric motive over $\Fp$, $i: \SpecFp \r \SpecZ$. More precisely, the special value prediction of \ref{conj_L} in this case is 
\eqn \mylabel{eqn_specialFp}
L^* (i_* N, 0) \equiv (\log p)^{- \chi(N\dual(-1))} \mod {\Q^\x},
\xeqn
where $\chi(N\dual(-1))$ is the Euler characteristic of motivic cohomology (see \refeq{Eulerchar}, computed in the category $\DMBeic(\Fp)$).
\xtheo
\pf
$\Rightarrow$: to show \refco{L} and \refeq{specialFp} for $i_* N$, we may replace $N$ by $\gr_*^W \pH^* N$, the weight graded pieces of the truncations with respect to the motivic $t$-structure, since both the weight filtration and the $t$-structure are bounded \refaxcohomdim. The subcategory of $\MM(\Fp)$ consisting of pure objects is, by \cite[Axiom 4.11]{Scholbach:fcoho}, the category of pure motives with respect to numerical equivalence, $\PureMot_\num(\Fp)$. Under \refco{numrat}, this agrees with Chow motives $\PureMot_\rat(\Fp)$. Finally, $\chi_{\DMBei(\Z)}((i_* N)\dual(-1)) = \chi_{\DMBei(\Fp)}(N\dual(-1))$, so we have to show $\ord_{s=0} L(i_* N) = - \dim \H^0(N\dual,1) = -\dim \H_0(N, -1)$ and $L^*(i_* N) \equiv (\log p)^{-\dim \H_0(N, -1)} \mod {\Q^\x}$.

Consider first $N = H := \M_{\Fp}(X)\twi{-n}$ with $X / \Fp$ smooth and projective. Then $L(i_* H) = L(\M_\Z(X)\twi{1-n})$. Let $Z^n(X) / \num$ be the group of codimension $n$ cycles on $X$ modulo numerical equivalence. Then 
$$\dim	 \H^0(H) = \rk \CH^n(X) \stackrel{\text{\ref{conj_numrat}}} = \rk Z^n(X) / \num = - \ord_{s=n} \zeta(X, s),$$ 
so the pole order claim holds for $H$ by assumption: the Tate conjecture and the agreement of the $\ell$-adic homological and numerical equivalence relations on $X$ (up to torsion) together are equivalent to the rightmost equality \cite[Thm. 2.9]{Tate:Conjectures}. 

In general, $N$ is a direct summand of $H$ as above. Let $N \oplus N' = H$, which as an object in $\PureMot_\rat(\Fp)$ is denoted $h(X)(n)$. By the previous case, 
\eqn \label{eqn_Tate}
\dim \H^0 N + \dim \H^0 N' = - \ord L(N) - \ord L(N'). 
\xeqn
Let $-_\ell: \PureMot_\rat(\Fp) \r \oplus \underline {\Ql[\Gal(\Fp)]}$, $\pi h(X)(n) \mapsto \oplus_a \pi^* \H^a(X, \Ql(n))$ be the $\ell$-adic realization functor taking values in graded continuous $\ell$-adic $\Gal(\Fp)$-representations. We write $\H^0(N_\ell) := N_\ell^{\Gal(\Fp)}$, the Galois cohomology of the $\ell$-adic Galois module $N_\ell$. 
The following way of reasoning is borrowed from \lcs We have the following chain of inequalities:
\eqnarr
- \ord_{s=0} L(N,s) & 
\geq & \dim_{\Ql} \ker (\Id - \Fr^{-1})|N_\ell \\
& \geq & \dim_{\Ql} (N_\ell)^{\Gal(\Fp)} \\
&=& \dim_{\Ql} \H^0 (N_\ell) \\
& \geq & \dim_\Q \H^0(N)
\xeqnarr
The last inequality is by the injectivity of the cycle class map $\H^0(N) \r \H^0(N_\ell)$, which follows from the injectivity of $\H^0(H) \r \H^0(H_\ell) = \H^{2n}(X, \Ql(n))$, i.e., the agreement of homological and rational equivalence, which holds under \refco{numrat}. Therefore, in \refeq{Tate} equality of dimensions must hold for the individual summands, so the pole order part is shown.

The claim \refeq{specialFp} and the special values formula of \ref{conj_L} trivially hold for $N = \one (-1)$: the residue of $L(i_* \one(-1), s) = \zeta(\SpecFp, s) = (1-p^{-s})^{-1}$ at $s=0$ is $(\log p)^{-1}$, which is the inverse of the determinant of $\pi_{\M(\Fp)} = \pi_{i_* i^* \one \twi{-1}}$ (\refex{Fp}). Jannsen's semisimplicity theorem for $\PureMot_\num(\Fp)$ yields a decomposition $N = \one(-1)^r \oplus R$ with $\Hom_{\PureMot_\num(\Fp)}(\one(-1), R) = \Hom_{\PureMot_\num(\Fp)}(R, \one(-1)) = 0$. Hence we can assume $N = R$. By the Lefschetz trace formula, the $L$-function of any pure motive over $\Fp$ is a rational function in $p^{-s}$ with rational coefficients that are independent of $\ell$, see e.g.\ \cite[Section 7.1.4]{Andre:Motifs}. By the preceding part, the $L$-function of $i_* R$ does not have a pole at $s = 0$, therefore the leading term of the Laurent series $L(i_* R, s)$ is simply the value at this point, a nonzero rational number (as opposed to an $\ell$-adic or, via $\sigma_\ell$,
 a complex number). 

$\Leftarrow$: we again use the theorem of Tate cited above: the Tate conjecture for $X / \Fp$ is implied by $\ord_{s=j} \zeta(X, s) = - \rk Z^j(X)/\num$. Under \ref{conj_numrat}, that term is $-\rk \CH^j(X) = -\dim \H^{2j}(\M(X)(-j))$. Thus, \refco{L} for $i_* \M(X)(-j)$ implies the Tate conjecture on the $j$-th Chow group of $X$.
\xpf

\bibliography{bib}  

\end{document}